\documentclass[11pt]{article}
\input amssym.def
\input amssym.tex
\usepackage{epsfig,psfig}
\newtheorem{theorem}{Theorem}

\def\binom#1#2{{#1}\choose{#2}}

\def\slashedfrac#1#2{\hbox{\kern.1em %
 \raise.5ex\hbox{\the\scriptfont0 #1}\kern-.11em %
 /\kern-.15em\lower.25ex\hbox{\the\scriptfont0 #2}}}

\newcommand{\eqn}[1]{(\ref{#1})}
\newcommand{\hsp}{\hspace*{\parindent}}

\newcommand{\eeq}{\end{equation}}
\newcommand{\beql}[1]{\begin{equation}\label{#1}}

\newcommand{\sO}{{\cal O}}
\newcommand{\sI}{{\cal I}}
\newcommand{\sW}{{\cal W}}

\newcommand{\sA}{{\cal A}}

\newcommand{\ZZ}{{\Bbb Z}}
\newcommand{\RR}{{\Bbb R}}

\newcommand{\CC}{{\Bbb C}}

\newcommand{\af}{\alpha}
\newcommand{\la}{\lambda}

\newcommand{\dd}{\ldots}

\newcommand{\sG}{{\cal G}}

\newcommand{\sD}{{\cal D}}

\newcommand{\sC}{{\cal C}}

\newcommand{\sP}{{\cal P}}
\newcommand{\sQ}{{\cal Q}}
\newcommand{\sT}{{\cal T}}

\makeatletter
\def\@sect#1#2#3#4#5#6[#7]#8{\ifnum #2>\c@secnumdepth
     \def\@svsec{}\else
     \refstepcounter{#1}\edef\@svsec{\csname the#1\endcsname.\hskip .75em }\fi
     \@tempskipa #5\relax
      \ifdim \@tempskipa>\z@
        \begingroup #6\relax
          \@hangfrom{\hskip #3\relax\@svsec}{\interlinepenalty \@M #8\par}%
        \endgroup
       \csname #1mark\endcsname{#7}\addcontentsline
         {toc}{#1}{\ifnum #2>\c@secnumdepth \else
                      \protect\numberline{\csname the#1\endcsname}\fi
                    #7}\else
        \def\@svsechd{#6\hskip #3\@svsec #8\csname #1mark\endcsname
                      {#7}\addcontentsline
                           {toc}{#1}{\ifnum #2>\c@secnumdepth \else
                             \protect\numberline{\csname the#1\endcsname}\fi
                       #7}}\fi
     \@xsect{#5}}
\def\@begintheorem#1#2{\it \trivlist \item[\hskip \labelsep{\bf #1\ #2.}]}

\def\plain{plain}\ifx\fmtname\plain\csname fi\endcsname
     
     \input docstrip
     \preamble

     Do not distribute the stripped version of this file.
     The checksum in the header refers to the documented version.

     \endpreamble
     \generateFile{here.sty}{t}{\from{here.doc}{}}
     \endinput
\fi
\ifcat a\noexpand @\let\next\relax\else\def\next{%
    \documentstyle[here,doc]{article}\MakePercentIgnore}\fi\next
\ifx\@Hxfloat\@Hundef\else\expandafter\endinput\fi
\let\@Hxfloat\@xfloat
\def\@xfloat#1[{\@ifnextchar{H}{\@HHfloat{#1}[}{\@Hxfloat{#1}[}}
\def\@HHfloat#1[H]{%
\expandafter\let\csname end#1\endcsname\end@Hfloat
\vskip\intextsep\vbox\bgroup\def\@captype{#1}\parindent\z@
\ignorespaces}
\def\end@Hfloat{\egroup\vskip \intextsep}

\makeatother

\catcode`\@=11
\renewcommand{\section}{
        \setcounter{equation}{0}
        \@startsection {section}{1}{\z@}{-3.5ex plus -1ex minus
        -.2ex}{2.3ex plus .2ex}{\large\bf}
        }
\catcode`@=12

\begin{document}
\begin{center}
{\Large {\bf Packing Lines, Planes, etc.: Packings in Grassmannian Spaces}} \\ [+.05in]
\vspace{1.5\baselineskip}
{\em J. H. Conway} \\
\vspace*{+.2\baselineskip}
Department of Mathematics \\
Princeton University, Princeton, NJ 08544 \\
\vspace*{+1\baselineskip}
{\em R. H. Hardin} and {\em N. J. A. Sloane} \\
\vspace*{+.2\baselineskip}
Mathematical Sciences Research Center \\
AT\&T Bell Laboratories \\
Murray Hill, NJ 07974 \\
\vspace{1\baselineskip}
April 12, 1996. Minor editorial changes, July 31, 2002. \\
\vspace{1.5\baselineskip}
{\bf ABSTRACT}
\vspace{.5\baselineskip}
\end{center}
\setlength{\baselineskip}{1.5\baselineskip}

This paper addresses the question:
how should $N$ $n$-dimensional subspaces of $m$-dimensional Euclidean
space be arranged so that they are as far apart as possible?
The results of extensive computations for modest values of $N, n,m$ are
described, as well as a reformulation of the problem that
was suggested by these computations.
The reformulation gives a way to describe $n$-dimensional subspaces of $m$-space as points on a sphere in dimension $(m-1) (m+2)/2$, which provides a (usually)
lower-dimensional representation than the Pl\"{u}cker embedding, and leads to
a proof that many of the new packings are optimal.
The results have applications to the graphical display of multi-dimensional
data via Asimov's ``Grand Tour'' method.
\clearpage
\thispagestyle{empty}
\setcounter{page}{1}

\section{Introduction}
\hsp
Although there is a considerable literature dealing with Grassmannian spaces
(see for example \cite{Cho49}, \cite{Lei61}, \cite{Won67},
\cite{JC74}, \cite{GH78}, \cite{GMM95}, \cite{Zan95}),
the problem of finding the best packings in such spaces seems to have
received little attention.

We have made extensive computations on this problem,
and have found a number of putatively optimal packings.
These computations have led us to conclude that the best definition of
distance on Grassmannian space is the ``chordal distance''
defined in Section~2.

Sections~3, 4, 5, 6 discuss the problems
of packing lines in $\RR^3$, planes in $\RR^4$, $n$-spaces in $\RR^m$, and
lines in $\RR^m$, respectively.

The complete results have been placed in a database on NJAS's home page.
Our search has concentrated on packings of $N\le 55$ subspaces of
$G(m,n)$, for $m \le 16$, $n \le 3$.

\section{Grassmannian space}
\hsp
The {\em Grassmannian space} $G(m,n)$ is the set of all $n$-dimensional
subspaces of real Euclidean $m$-dimensional space $\RR^m$.
This is a homogeneous space isomorphic to $O(m) / (O(n) \times O(m-n))$, and forms a compact Riemannian manifold of dimension $n(m-n)$.

We first discuss how to define the distance between two $n$-planes
$P$, $Q \in G(m,n)$.
The {\em principal angles} $\theta_1, \ldots , \theta_n \in [0, \pi/2 ]$ between
$P$ and $Q$ are defined by (we follow \cite{GVL89}, p.~584)
$$\cos \theta_i = \max_{u \in P} \max_{v \in Q} u \cdot v = u_i \cdot v_i ~,$$
for $i=1, \ldots , n$,
subject to $u \cdot u = v \cdot v =1$, $u \cdot u_j =0$, $v \cdot v_j =0$ $(1 \le j \le i-1)$.
The vectors $\{u_i\}$ and $\{v_j \}$ are {\em principal vectors}
corresponding to the pair $P$ and $Q$.

\begin{figure}[htb]
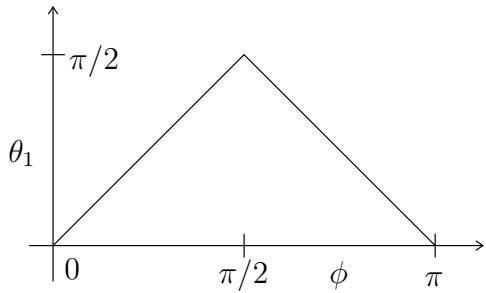

\begin{center}
\input tent_text.tex
\end{center}
\caption{Principal angle $\theta_1$ between two lines as the angle between them increases from 0 to $\pi$.}
\label{fg1}
\end{figure}

Wong \cite{Won67} shows that the {\em geodesic distance} on $G(m,n)$ between $P$ and $Q$ is
\beql{Eq1}
d_g (P,Q) = \sqrt{\theta_1^2 + \cdots + \theta_n^2} ~.
\eeq
However, this definition has one drawback:
it is not everywhere differentiable.
Consider the case $n=1$, for example, and hold one line $P$ fixed while rotating another line $Q$ (both passing through the origin).
As the angle $\phi$ between $P$ and $Q$ increases from 0 to $\pi$,
the principal angle $\theta_1$
increases from 0 to $\pi /2$ and then falls to 0, and is non-differentiable at $\pi /2$ (see Fig.~\ref{fg1}).

Although one might expect this non-differentiability to be a mere technicality, it does
in fact cause considerable difficulties for our optimizer, especially
in higher dimensions in cases when many distances fall in the neighborhood of singular
points of $d_g$.

An alternative measure of distance, which we call the {\em chordal distance}, is given by
\beql{Eq2}
d_c (P,Q) = \sqrt{\sin^2 \theta_1 + \cdots + \sin^2 \theta_n} ~.
\eeq
The reason for the name will be revealed later.
This approximates the geodesic distance when the planes are close,
has the property that its square is differentiable everywhere,
and, as we shall attempt to demonstrate,
has a number of other desirable features.

A third definition has been used by Asimov \cite{Asi85} and Golub and Van~Loan \cite{GVL89}, p.~584,
namely
$$d_m (P,Q) = \theta_n = \max_{i=1, \dd, n} \theta_i ~.$$
This shares the vices of the geodesic distance.

Of course for $n=1$ all three definitions are equivalent, in the sense that
they lead to the same optimal packings.

We can now state the packing problem:
given $N,n,m$, find a set of $n$-planes $P_1, \dd, P_N \in G(m,n)$ so that
$\min\limits_{i \neq j} d(P_i, P_j)$ is as large as possible,
where $d$ is either geodesic or chordal distance.
Since $G(m,n)$ is compact, the problem is well-defined.
Because $G(m,n)$ and $G(m, m-n)$ are essentially the same space, we will usually assume that $n \le m/2$.

We also need some further terminology.
A {\em generator matrix} for an $n$-plane $P\in G (m,n)$ is an $n \times m$ matrix whose rows span $P$.
The orthogonal group $O(m)$ acts on $G(m,n)$ by right multiplication of generator matrices.
The {\em automorphism group} of a subset $\{P_1, \dd, P_N \} \subset G(m,n)$ is the subset of $O(m)$ which fixes or permutes these planes.

By applying a suitable element of $O(m)$ and choosing appropriate
basis vectors for the planes, any given pair of $n$-planes $P,Q$ with $n \le m/2$ can be assumed to have generator matrices
\beql{Eq3}
\left[
\matrix{
1 & 0 & \cdots & 0 & 0 & 0 & \cdots & 0 & 0 & \cdots & 0 \cr
0 & 1 & \cdots & 0 & 0 & 0 & \cdots & 0 & 0 & \cdots & 0 \cr
\cdot & \cdot & \cdots & \cdot & \cdot & \cdot & \cdots & \cdot & \cdot & \cdots & \cdot \cr
0 & 0 & \cdots & 1 & 0 & 0 & \cdots & 0 & 0 & \cdots & 0\cr}
\right]
\eeq
and
\beql{Eq4}
\left[
\matrix{\cos \theta_1 & 0 & \cdots & 0 & \sin \theta_1 & 0 & \cdots & 0 & 0 & \cdots & 0 \cr
0 & \cos \theta_2 & \cdots & 0 & 0 & \sin \theta_2 & \cdots & 0 & 0 & \cdots & 0 \cr
\cdot & \cdot & \cdots & \cdot & \cdot & \cdot & \cdots & \cdot & \cdot & \cdots & \cdot \cr0 & 0 & \cdots & \cos \theta_n & 0 & 0 & \cdots & \sin \theta_n & 0 & \cdots & 0 \cr}
\right]
\eeq
respectively, where $\theta_1 , \dd, \theta_n$ are the principal angles between them
(\cite{Won67}, Theorem~2).
\section{Packing lines through the origin in 3-space}
\hsp
Our initial work on this problem was prompted by a question raised in 1992 by Julian Rosenman, an oncologist at the Univ. of North Carolina School of Medicine and Computer Science,
in connection with the treatment of tumors using high energy laser beams \cite{Ros92}.
He asked for the best way to separate $N$ lines through a given point in $\RR^3$, or in other words for the best packings in $G(3,1)$.

Together with W.~D. Smith, we had been carrying out an extensive search for the best
packings of a given number of points on $S^2$, i.e. {\em spherical codes} \cite{HSS94},
\cite{HSS96}, and we therefore
modified our programs to search instead for packings of lines.
We omit the details of this search, since we later
recalculated these results using the more general methods described in Section~4.

The results are summarized in Table~\ref{ta1}, which gives, for each value of $N$
in the range 2 to 55, the minimal angle of the best packing we have found of $N$ lines,
and for comparison the minimal angle of the best packing known of $2N$ points
on $S^2$ (taken from \cite{HSS94}).
We see that requiring a packing of $2N$ points on $S^2$ to be antipodal is a definite handicap:
only in the cases of 6 and 12 points do the antipodal and unrestricted packings coincide.
Decimals in the tables have been rounded to four decimal places.

The last two columns of the table specify the largest automorphism group\footnote{That is, the subgroup of $O(3)$ that fixes or permutes the $2N$ points.} we have found of any such best antipodal packing of $2N$ points.
The fourth column gives the order of the group and the fifth and sixth columns its name in the orbifold notation
(cf. \cite{CS96}) and as the double cover of a rotation group.
The symbol $\pm \sG$ indicates that the group consists of the matrices
$\pm M$ for $M \in \sG$, where $\sG$ is a cyclic $( \sC )$, dihedral $( \sD )$,
tetrahedral $(\sT)$,
octahedral $(\sO)$ or icosahedral $(\sI)$ group.
In each case the subscript gives the order of the rotation group.
\begin{table}[H]
\caption{Best packings found of $N$ lines through origin in $\RR^3$.
(The third column gives the best packing known of $2N$ points on a sphere.)}

$$
\begin{array}{rlrrlll}
\multicolumn{1}{c}{\mbox{No. of}} & 
\multicolumn{1}{c}{\mbox{Min.}} &
\multicolumn{1}{c}{\mbox{Min.}} &
\multicolumn{1}{c}{\mbox{Group}} &
\multicolumn{1}{c}{\mbox{Group}} &
\multicolumn{1}{c}{\mbox{Group}} &
\multicolumn{1}{c}{\mbox{Notes}} \\
\multicolumn{1}{c}{\mbox{lines}} &
\multicolumn{1}{c}{\mbox{angle}} &
\multicolumn{1}{c}{\mbox{angle}} &
\multicolumn{1}{c}{\mbox{order}} &
\multicolumn{1}{c}{\mbox{name}} & \multicolumn{1}{c}{\mbox{structure}} \\
\multicolumn{1}{c}{N} & ~ &
\multicolumn{1}{c}{\mbox{(packing)}} \\ \hline
2~~~ & 90.0000 & 109.4712 & 16 & *224 & \pm \sD_8 & \mbox{square} \\
3~~~ & 90.0000 & 90.0000 & 48 & *432 & \pm \sO_{24} & \mbox{octahedron} \\
4~~~ & 70.5288 & 74.8585 & 48 & *432 & \pm \sO_{24} & \mbox{cube} \\
5~~~ & 63.4349 & 66.1468 & 20 & 2\! * \!5 & \pm \sD_{10} & \mbox{pentagonal antiprism (see note)} \\
6~~~ & 63.4349 & 63.4349 & 120 & *532 & \pm \sI_{60} & \mbox{icosahedron} \\
7~~~ & 54.7356 & 55.6706 & 48 & *432 & \pm \sO_{24} & \mbox{rhombic dodecahedron} \\
8~~~ & 49.6399 & 52.2444 & 2 & \times & \pm \sC_1 & \mbox{see note} \\
9~~~ & 47.9821 & 49.5567 & 6 & 3\times & \pm \sC_3 & \mbox{see note} \\
9~~~ & 47.9821 & 49.5567 & 4 & 2* & \pm \sC_2 & ~ \\ 
10~~~ & 46.6746 & 47.4310 & 24 & *226 & \pm \sD_{12} & \mbox{hexakis bi-antiprism (see note)}\\
11~~~ & 44.4031 & 44.7402 & 20 & 2\! * \! 5 & \pm \sD_{10} & ~ \\
12~~~ & 41.8820 & 43.6908 & 48 & *432 & \pm \sO_{24} & \mbox{rhombicuboctahedron} \\
13~~~ & 39.8131 & 41.0377 & 4 & 2* & \pm \sC_2 & ~ \\
14~~~ & 38.6824 & 39.3551 & 2 & \times & \pm \sC_1 & ~  \\
15~~~ & 38.1349 & 38.5971 & 20 & 2 \!* \!5 & \pm \sD_{10} & \mbox{see note} \\
16~~~ & 37.3774 & 37.4752 & 120 & *532 & \pm \sI_{60} & \mbox{pentakis dodecahedron} \\ 
17~~~ & 35.2353 & 35.8078 & 2 & \times & \pm \sC_1 &  \\
18~~~ & 34.4088 & 35.1897 & 6 & 3\times & \pm \sC_3 &  \\
19~~~ & 33.2115 & 34.2507 & 2 & \times & \pm \sC_1 &  \\
20~~~ & 32.7071 & 33.1584 & 8 & *222 & \pm \sD_4 &  \\
21~~~ & 32.2161 & 32.5064 & 10 & 5\times & \pm \sC_5  \\
22~~~ & 31.8963 & 31.9834 & 12 & 2\! * \!3 & \pm \sD_6  \\
23~~~ & 30.5062 & 30.9592 & 2 & \times & \pm \sC_1 &  \\
24~~~ & 30.1628 & 30.7628 & 24 & 3\! * \!2 & \pm \sT_{12}  \\
25~~~ & 29.2486 & 29.7530 & 6 & 3\times & \pm \sC_3 & \\
26~~~ & 28.7126 & 29.1948 & 4 & 2* & \pm \sC_2 & \\
27~~~ & 28.2495 & 28.7169 & 2 & \times & \pm \sC_1  \\
\end{array}
$$
\label{ta1}
\end{table}
\begin{table}[H]
\begin{center}
Table~1 (cont.): Best packings found of $N$ lines through origin in $\RR^3$. \\
~ \\
\end{center}

$$
\begin{array}{rlrrlll}
\multicolumn{1}{c}{\mbox{No. of}} &
\multicolumn{1}{c}{\mbox{Min.}} &
\multicolumn{1}{c}{\mbox{Min.}} &
\multicolumn{1}{c}{\mbox{Group}} &
\multicolumn{1}{c}{\mbox{Group}} &
\multicolumn{1}{c}{\mbox{Group}} &
\multicolumn{1}{c}{\mbox{Notes}} \\
\multicolumn{1}{c}{\mbox{lines}} &
\multicolumn{1}{c}{\mbox{angle}} &
\multicolumn{1}{c}{\mbox{angle}} &
\multicolumn{1}{c}{\mbox{order}} &
\multicolumn{1}{c}{\mbox{name}} & \multicolumn{1}{c}{\mbox{structure}} \\
\multicolumn{1}{c}{N} & ~ &
\multicolumn{1}{c}{\mbox{(packing)}} \\ \hline
28~~~ & 27.8473 & 28.1480 & 2 & \times & \pm \sC_1  \\
29~~~ & 27.5244 & 27.5564 & 28 & 2\!*\!7 & \pm \sD_{14}  \\
30~~~ & 26.9983 & 27.1928 & 20 & 2\! * \!5 & \pm \sD_{10}  \\
31~~~ & 26.4987 & 26.6840 & 10 & 5\times & \pm \sC_5  \\
32~~~ & 25.9497 & 26.2350 & 4 & 2* & \pm \sC_2  \\
33~~~ & 25.5748 & 25.9474 & 24 & 3\! * \!2 & \pm \sT_{12}  \\
34~~~ & 25.2567 & 25.4638 & 6 & 3\times & \pm \sC_3  \\
35~~~ & 24.8702 & 25.1709 & 2 & \times & \pm \sC_1  \\
36~~~ & 24.5758 & 24.9265 & 6 & 3\times & \pm \sC_3  \\
37~~~ & 24.2859 & 24.4209 & 2 & \times & \pm \sC_1  \\
38~~~ & 24.0886 & 24.1282 & 2 & \times & \pm \sC_1  \\
39~~~ & 23.8433 & 23.9310 & 6 & 3\times & \pm \sC_3  \\
40~~~ & 23.3293 & 23.5531 & 2 & \times & \pm \sC_1  \\
41~~~ & 22.9915 & 23.1946 & 2 & \times & \pm \sC_1  \\
42~~~ & 22.7075 & 23.0517 & 6 & 3\times & \pm \sC_3  \\
43~~~ & 22.5383 & 22.6744 & 12 & 2\! * \!3 & \pm \sD_6  \\
44~~~ & 22.2012 & 22.4679 & 2 & \times & \pm \sC_1  \\
45~~~ & 22.0481 & 22.1540 & 4 & 2* & \pm \sC_{2}  \\
46~~~ & 21.8426 & 22.0276 & 2 & \times & \pm \sC_1  \\
47~~~ & 21.6609 & 21.7221 & 2 & \times & \pm \sC_1  \\
48~~~ & 21.4663 & 21.5206 & 6 & 3\times & \pm \sC_3  \\
49~~~ & 21.1610 & 21.3711 & 2 & \times & \pm \sC_1  \\
50~~~ & 20.8922 & 21.0312 & 2 & \times & \pm \sC_1  \\
51~~~ & 20.6903 & 20.8556 & 2 & \times & \pm \sC_1  \\
52~~~ & 20.4914 & 20.6566 & 2 & \times & \pm \sC_1  \\
53~~~ & 20.2685 & 20.4394 & 2 & \times & \pm \sC_1  \\
54~~~ & 20.1555 & 20.3044 & 6 & 3\times & \pm \sC_3  \\
55~~~ & 20.1034 & 20.1110 & 120 & *532 & \pm \sI_{60} & \mbox{see note}  \\
\end{array}
$$
\end{table}

In some cases the best packings can be obtained by taking the diameters of a known
polyhedron, and if so this is indicated in the final column of the table.
Other entries in this column refer to the brief descriptions given below.

The entries for $N \le 6$ were shown to be optimal by Fejes T\'{o}th in 1965 \cite{FT65}
(see also Rosenfeld \cite{Ros94}),
and the 7-line arrangement will be proved optimal in Section~5.
The solutions for $N \ge 8$ are the best found with over 15000
random starts with our optimizer.
There is no guarantee that these are optimal, but experience with similar problems
suggests that they will be hard to beat and in any case will be not far from optimal.

For $N=1,2,3,6$ the solutions are known to be unique,
for $N=4$ there are precisely two solutions (\cite{FT65}, \cite{Ros93}, \cite{Ros94}), and for $N=5,7,8$ the solutions appear to be unique.
For larger values of $N$, however, the solutions are often not unique.
For $N=9$ lines there are two different solutions, and in the range $N \le 30$ the solutions
for 10, 22, 25, 27, 29 lines (and possibly others)
contain lines that ``rattle'', that is, lines which can be moved freely over a small range of angles
without affecting the minimal angle.
\subsection*{Notes on Table~1}
\begin{figure}[htb]
\centerline{\psfig{file=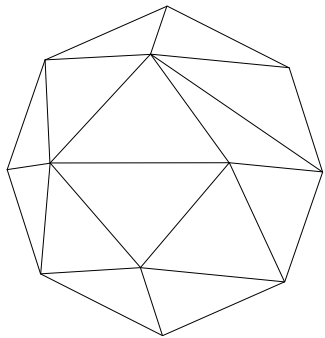,width=1.9in}}
\caption{Best antipodal packing found of 16 points.}
\label{fg2}
\end{figure}

\begin{figure}[htb]
\centerline{\psfig{file=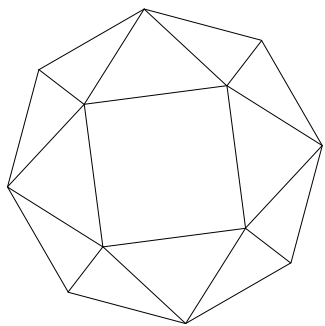,width=1.9in}}
\caption{Best (unrestricted) packing of 16 points known.}
\label{fg3}
\end{figure}

\begin{itemize}
\item[$N=5$.]
Five of the six diameters of a regular icosahedron.
\item[$N=8$.]
The putatively optimal arrangement forms an unpleasant-looking
configuration of 16 antipodal points on $S^2$, with no further
symmetry.
The convex hull is shown in Fig.~\ref{fg2}.
In contrast, the putatively best packing of 16 points (Fig.~\ref{fg3}) has a group of order 16.
\end{itemize}

\begin{figure}[htb]
\centerline{\psfig{file=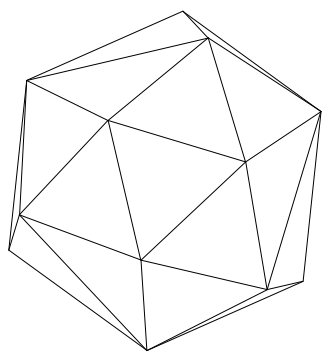,width=1.9in}}
\caption{Antipodal packing of 18 points with group of order 6.}
\label{fg4}
\end{figure}

\begin{figure}[htb]
\centerline{\psfig{file=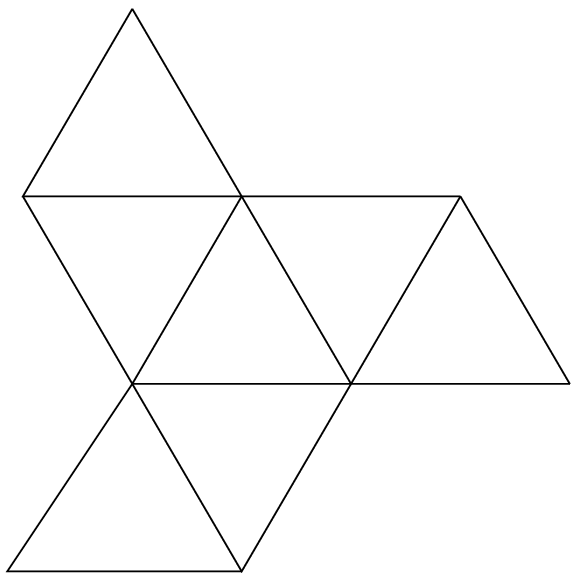,width=1.9in}}
\caption{}
\label{fg5}
\end{figure}

\begin{figure}[htb]
\centerline{\psfig{file=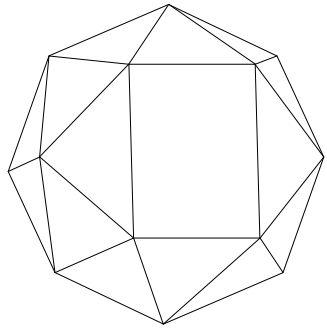,width=1.9in}}
\caption{Antipodal packing of 18 points with group of order 4.}
\label{fg6}
\end{figure}

\begin{figure}[htb]
\centerline{\psfig{file=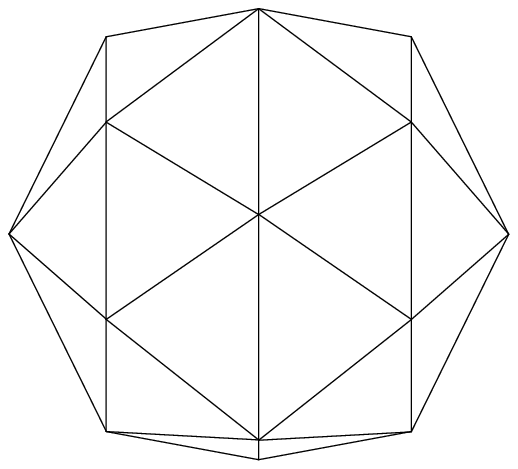,width=1.9in}}
\caption{Best antipodal packing found of 20 points.}
\label{fg7}
\end{figure}

\begin{figure}[htb]
\centerline{\psfig{file=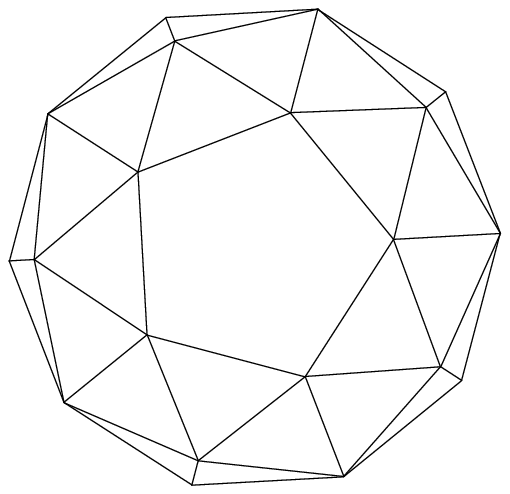,width=1.9in}}
\caption{Best antipodal packing found of 30 points.}
\label{fg8}
\end{figure}

\begin{figure}[htb]
\centerline{\psfig{file=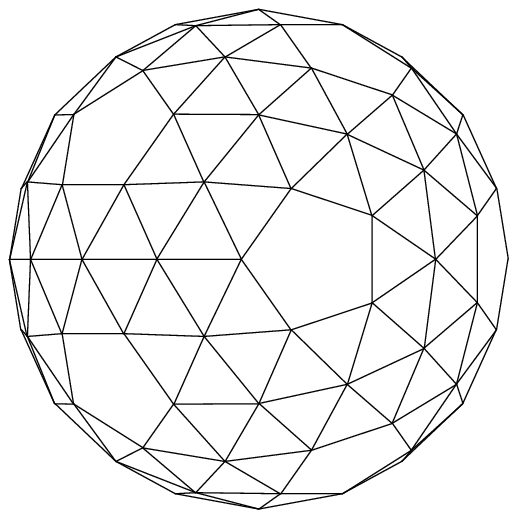,width=1.9in}}
\caption{Best antipodal packing found of 110 points.}
\label{fg9}
\end{figure}

\paragraph{$N=9$.}
There appear to be two inequivalent solutions.
The nicest has symmetry group $3 \times$, of order 12 (or $[2^+,6^+]$ in
Coxeter's notation \cite{CM80}) --- see Fig.~\ref{fg4}.
The points lie in six horizontal layers of equilateral triangles.
The nine points in the Northern hemisphere
can be located by drawing seven equilateral spherical triangles
of edge length 47.98213264$^\circ$, arranged as in Fig.~\ref{fg5}.
The North pole is located at the midpoint of the central triangle.
The second solution, shown in Fig.~\ref{fg6}, has group $2*$, of order 4.
\paragraph{$N=10$.}
This is a ``hexakis bi-antiprism'', since it consists of a bi-(hexagonal antiprism)
together with an axial line that can ``rattle'',
giving infinitely many solutions.
The axis appears horizontally in Fig.~\ref{fg7}.
\paragraph{$N=15$.}
Combinatorially, this is an ``axially depleted pentakis
dodecahedron'',
obtained from a pentakis dodecahedron (the solution
for $N=16$) by omitting two opposite vertices (see Fig.~\ref{fg8}).
However, the angular separation is slightly greater than for $N=16$.
\paragraph{$N=55$.}
The 110 antipodal points can be taken to be the union of the vertex sets of
a dodecahedron (20), an icosidodecahedron (30) and a truncated
icosahedron (60) --- see Fig.~\ref{fg9}.
The 15 lines through the icosidodecahedral points ``rattle''.
This arrangement of points can be obtained from a ``geodesic dome'' or
12-fold reticulated icosahedron \cite{Cox74} by omitting its twelve pentagonal vertices,
and so may be described as consisting of the vertices of a ``depleted 12-fold reticulated icosahedron''.
Incidentally, our best solution for 61 lines has a symmetry group
of order 6, and is not found by placing points at the centers of the pentagonal faces of
Fig.~\ref{fg9}.

Packings of lines in higher dimensional spaces will be discussed in Section~6.

\section{Packing two-dimensional planes in $\RR^4$}
\hsp
In \cite{Asi85} Asimov proposed a technique called the ``Grand Tour'' for
displaying multi-dimensional data on a two-dimensional computer screen.
His idea was to choose a finite sequence or ``tour'' of two-dimensional planes that are in some sense ``representative'' of all planes, and to project
the data onto each plane in turn,
in the hope that the viewer will be able to notice any pattern
or structure that is present.
This technique has been implemented in the XGobi program
\cite{BS86}, \cite{SCB91}.

In 1993 Dianne Cook asked us if we could modify our algorithm for finding
spherical codes, in order to search for packings in $G(m,2)$ for $m \ge 4$.
Furthermore, for the Grand Tour application, there should be a Hamiltonian
circuit through the planes, in which edges
indicate neighboring planes.
We accomplished this in the following way.

Following the methods we had used to find packings in
$S^d$ and experimental designs in various spaces (see \cite{HS92},
Sec.~4.11;
\cite{HS93}, \cite{HSS96}), corresponding to a set of planes $S= \{P_1, \dd, P_N \} \subseteq G(m,n)$ we define a potential
$$\Phi_c (S) = \sum_{i<j}
\frac{1}{d_c (P_i, P_j) -A} ~,$$
where $A$ is a suitably chosen constant.
There is a similar definition for $\Phi_g$ involving $d_g$.
Initially $A$ is set to 0 and $S$ is a randomly chosen set of planes.
We invoke the Hooke-Jeeves ``pattern search'' optimization
algorithm to modify $S$, attempting to minimize $\Phi_c (S)$.
After a fixed time (we used 100 steps of the optimizer),
$A$ is advanced halfway to the current minimal distance between the $P_i$,
and the process is repeated, terminating when no further improvement
is obtained to the accuracy of the machine.
The whole process is repeated with several thousand random starts --- and also with initial configurations taken from other sources, such as subsets or supersets of other arrangements --- and the best final configuration is recorded.

The partial derivatives of $\Phi_c (S)$ were found analytically, but for
$\Phi_g (S)$ we calculated them by numerical differentiation.

We began by computing a table of packings of $N$ planes in $G(4,2)$, for $N \le 50$, with respect
to both $d_c$ and $d_g$.
The results for packings of $N \le 24$ planes are summarized in Table~\ref{ta2}.
Columns 2 and 3 refer to the best packings found for the chordal distance, and give both
$$d_c^2 = \min_{i \neq j} d_c^2 (P_i, P_j) ~~{\rm and}~~
d_g^2 = \min_{i \neq j} d_g^2 (P_i, P_j ) $$
for such packings.
Columns 4 and 5 refer similarly to the best packings found with respect
to geodesic distance.
\begin{table}[htb]
\caption{Comparison of putatively best packings of $N$ planes in $G(4,2)$ with respect to $d_c$ and $d_g$.}

$$
\begin{array}{rcclcc}
~ & \multicolumn{2}{c}{\mbox{Best w.r.t. $d_c$}} & ~~~~ &
\multicolumn{2}{c}{\mbox{Best w.r.t. $d_g$}} \\ [+.1in]
N & d_c^2 & d_g ^2 & ~~~~~~~ & d_c^2 & d_g ^2 \\ [+.1in]
2 & 2.0000 & 4.9348 &  & 2.0000 & 4.9348 \\
3 & 1.5000 & 2.1932 &  & 1.2500 & 2.7416 \\
4 & 1.3333 & 1.8253 &  & 1.2000 & 2.6824 \\
5 & 1.2500 & 1.9739 &  & 1.2000 & 2.6824 \\
6 & 1.2000 & 2.6824 &  & 1.2000 & 2.6824 \\
7 & 1.1667 & 1.6440 &  & 0.9875 & 2.1281 \\
8 & 1.1429 & 1.5818 &  & 0.9700 & 1.9235 \\
9 & 1.1231 & 1.5175 &  & 0.9501 & 1.8087 \\
10 & 1.1111 & 1.5725 &  & 0.9764 & 1.7886 \\
11 & 1.0000 & 1.2715 &  & 0.9247 & 1.6711 \\
12 & 1.0000 & 1.3413 &  & 0.9204 & 1.6416 \\
13 & 1.0000 & 1.2348 &  & 0.9133 & 1.5600 \\
14 & 1.0000 & 1.2337 &  & 0.8933 & 1.5327 \\
15 & 1.0000 & 1.2337 &  & 0.8923 & 1.5284 \\
16 & 1.0000 & 1.2337 &  & 0.8904 & 1.5210 \\
17 & 1.0000 & 1.2337 &  & 0.8549 & 1.3925 \\
18 & 1.0000 & 1.2337 &  & 0.8504 & 1.3768 \\
19 & 0.9091 & 1.1666 &  & 0.8412 & 1.3477 \\
20 & 0.9091 & 1.1666 &  & 0.8351 & 1.3284 \\
21 & 0.8684 & 1.0352 &  & 0.8225 & 1.2834 \\
22 & 0.8629 & 1.0592 &  & 0.8046 & 1.2385 \\
23 & 0.8451 & 1.0081 &  & 0.7910 & 1.2012 \\
24 & 0.8372 & 0.9901 &  & 0.7812 & 1.1707
\end{array}
$$
\label{ta2}
\end{table}

Understanding these results was hindered by the fact that planes
are, as the name suggests, {\em plain},
with no distinguishing features, and when produced by
the computer appear as random generator matrices referred to a random coordinate frame.
However, we found that the set of principal vectors in the planes could often be used
to find a coordinate system that would reveal some of the structure of the planes.

For two planes the best packing is the same for both metrics, and consists of
two mutually orthogonal planes with principal angles $\pi/2$, $\pi/2$,
so $d_c^2 =2$, $d_g^2 = \pi^2/2$.

For three planes the two answers are different.
The best packing for chordal distance (Fig.~\ref{fg10}a) has principal
angles $\pi /3$, $\pi /3$ between each pair of planes,
so $d_c^2 = 3/2$, $d_g^2 = 2 \pi^2 /9$, whereas the best packing
for geodesic distance (Fig.~\ref{fg10}b) has angles $\pi /6$, $\pi /2$,
so $d_c^2 = 5/4$, $d_g^2 = 5 \pi^2 /18$.
\begin{figure}[htb]
$$
\begin{array}{cccc} \hline
1 & 0 & 1 & 0 \\
0 & 1 & 0 & 1 \\ \hline
1 & 0 & - \frac{1}{2} & \frac{\sqrt{3}}{2} \\
0 & 1 & - \frac{\sqrt{3}}{2} & - \frac{1}{2} \\ \hline
1 & 0 & - \frac{1}{2} & - \frac{\sqrt{3}}{2} \\
0 & 1 & \frac{\sqrt{3}}{2} & - \frac{1}{2} \\ \hline
\multicolumn{4}{c}{~~~} \\
\multicolumn{4}{c}{{\rm (a)}}
\end{array}
~~~~~~~~~~~~
\begin{array}{cccc} \hline
1~~ & 0 & 0 & 0 \\
0~~ & - \frac{1}{\sqrt{2}} & 0 & \frac{1}{\sqrt{2}} \\ \hline
0~~ & 1 & 0 & 0 \\
0~~ & 0 & - \frac{1}{\sqrt{2}} & \frac{1}{\sqrt{2}} \\ \hline
0~~ & 0 & 1 & 0 \\
- \frac{1}{\sqrt{2}} & 0 & 0 & \frac{1}{\sqrt{2}} \\ \hline
\multicolumn{4}{c}{~~~} \\
\multicolumn{4}{c}{{\rm (b)}}
\end{array}
$$
\caption{Generator matrices for best packings of three planes
in $G(4,2)$ with respect to (a) chordal and (b) geodesic distance.}
\label{fg10}
\end{figure}

Postponing discussion of $N=4$ and 5 for the moment, let us
consider the case of
six planes,
where we discovered that the answer for both metrics formed a regular simplex\footnote{Although not an isoclinic configuration, cf. Wong \cite{Wong61}.
There seems to be only a slight overlap connection between that problem and ours.}:
the principal angles between every pair of planes were arcsin $1/ \sqrt{5}$ and $\pi /2$.
This set of planes is conveniently described using
simplicial coordinates.
Let $A$, $B$, $C$, $D$, $E$ be the vectors from the center of a
regular simplex in $\RR^4$ to its vertices,
with $A+B+C+D+E =0$, and write
$[a,b,c,d,e]$ for $aA + bB + cC + dD + eE$.
Then one of the six planes is spanned by $[1, \tau , 1, 0,0]$ and its cyclic shifts, where
$\tau = (1+ \sqrt{5} )/2$, and
the other five planes are obtained from
it by taking even permutations of these coordinates.
(We later found an equally good packing with respect to chordal distance that was not a simplex with
respect to geodesic distance.
This is described below.)

The six planes intersect the surface of the unit ball of $\RR^4$
in a remarkable link, the
\linebreak
``Hexalink''.
It can be shown that apart from the $n$-component unlinks and Hopf links, which exist for all $n$,
the Hexalink is the only link of circular rings in which there are orientation-preserving symmetries taking any two links
to any other two.
(Details will be given elsewhere.)

The discovery of this 6-vertex simplex was initially somewhat of a surprise,
since $G(4,2)$ is only a 4-dimensional
manifold.
It did suggest that the answer should somehow be related to the six
equi-angular diameters of the icosahedron, and led to the following reformulation of the problem.

We remind the reader that any element $\alpha$ of $SO(4)$ may be represented
as
$$\alpha : x \mapsto \bar{\ell} xr ~,$$
where $x=x_0 + x_1 i + x_2 j + x_3 k$ represents a point on $S^3$ and
$\ell$, $r$ are unit quaternions \cite{DV64}.
The pair $- \ell$, $-r$ represent the same $\alpha$.
The correspondence between $\alpha$ and $\pm ( \ell , r)$ is one-to-one.

Given a plane $P \in G(4,2)$, let $\alpha$ be the element of $SO(4)$ that fixes $P$ and negates the points of the orthogonal plane $P^\perp$.
Then $\alpha^2 =1$, and for this $\alpha$,
it is easy to see that
$\ell = \ell_1 i + \ell_2 j + \ell_3 k$ and $r= r_1 i + r_2 j + r_3 k$ are purely
imaginary unit quaternions.
This establishes the following result (which can be found for example
in Leichtweiss \cite{Lei61}).
\begin{theorem}\label{th1}
A plane $P \in G (4,2)$ is represented by a pair $(\ell,r) \in S^2 \times S^2$, with $(- \ell , -r)$ representing the same $P$.
\end{theorem}

There are simple formulae relating $P$ and $(\ell ,r)$,
pointed out to us by Simon Kochen.
Given $(\ell,r)$, if $\ell \neq -r$ then $P$ is spanned by the vectors
corresponding to the quaternions $u=1 - \ell r$ and $v= \ell + r$.
The special case when $\ell =r$ is even simpler:
take $u=1$, $v= \ell$.
If $\ell = -r$, take $u$ and $v$ to be purely imaginary unit quaternions such that $\ell , u, v$ correspond to a coordinate frame.
Conversely, if $P$ is spanned by two orthogonal unit vectors
represented by quaternions $u$, $v$, then $\pm ( \ell , r) = \pm (u \bar{v} - v \bar{u} , \bar{v} u - \bar{u} v )$.

Given two planes $P,Q \in G(4,2)$, represented by $\pm ( \ell , r)$, $\pm ( \ell ' , r' )$, respectively, the principal angles $\theta_1$, $\theta_2$ between them may be found as follows.
Let $\phi$ (resp. $\psi$) be the angle between $\ell$ and $\ell '$ (resp. $r$ and $r'$), with $0 \le \phi , \psi \le \pi$.
If $\phi + \psi > \pi$, replace $\phi$ by $\pi - \phi$ and $\psi$ by
$\pi - \psi$,
so that $0 \le \phi + \psi \le \pi$, with $\phi \le \psi$ (say).
Then
\beql{Eq10}
\theta_1 , \theta_2 = \frac{\psi \pm \phi}{2} ~,
\eeq
\begin{eqnarray}
\label{Eq11}
d_g^2 (P,Q) & = & \frac{\psi^2 + \phi^2}{2} ~, \\
\label{Eq12}
d_c^2 (P,Q) & = & 1- \cos \psi ~ \cos \phi ~.
\end{eqnarray}
We omit the elementary proof.

A set $S = \{ P_1, \dd, P_N \} \subseteq G(4,2)$ is thus represented by
a ``binocular code'' consisting of a set of pairs
$\pm ( \ell_i , r_i) \in S^2 \times S^2$.
We call the list of $2N$ points $\pm \ell_i$ (they need not be distinct)
the ``left code'' corresponding to $S$, and the points $\pm r_i$ the ``right code''.
Conversely, given two multisets $L \subseteq S^2$, $R \subseteq S^2$,
each of size $2N$ and
closed under negation, and a bijection or ``matching'' $f$ between
them that satisfies $f( - \ell ) = -f ( \ell )$, $\ell \in L$, we obtain a set of $N$ planes in $G(4,2)$.

\begin{figure}[htb]
\centerline{\psfig{file=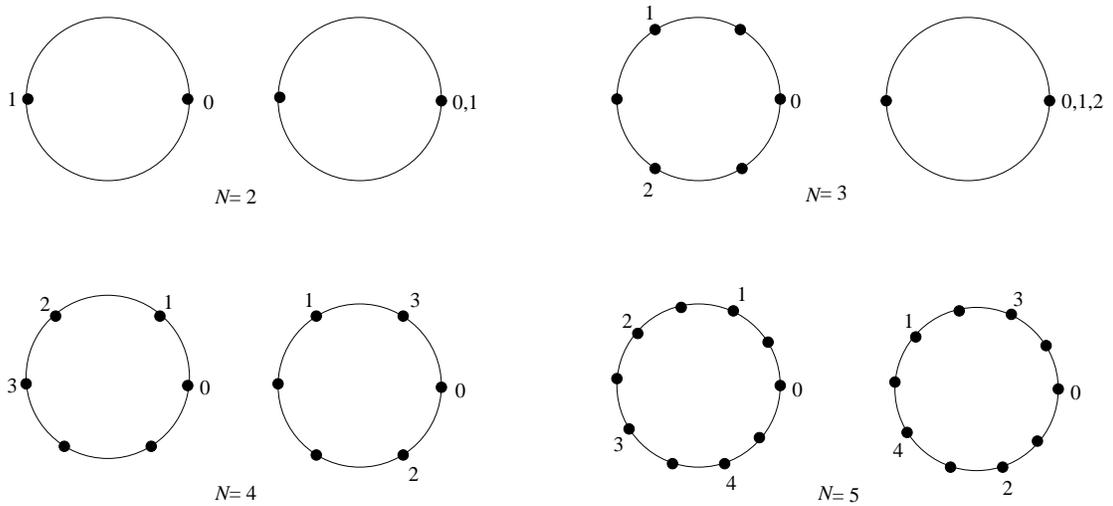,width=5.75in}}

\caption{Binocular codes describing best packings of $N=2 , \dd , 5$ planes in $G(4,2)$ for chordal distance.}
\label{fg11}
\end{figure}

\begin{figure}[htb]
\centerline{\psfig{file=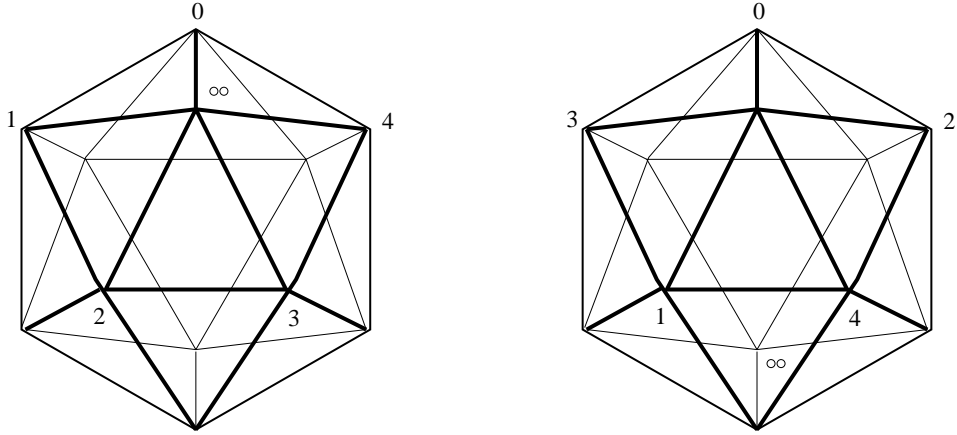,width=5in}}

\caption{Best packing of 6 planes in $G(4,2)$ with respect to both
metrics.
Left and right codes comprise vertices of icosahedron.
Adjacent vertices in one code are matched with non-adjacent vertices in the other code.}
\label{fg12}
\end{figure}

The binocular codes for the $d_c$-optimal packings of $N=2 , \dd , 6$ planes are shown in Figs.~\ref{fg11}, \ref{fg12}.
Except for $N=3$, the left and right codes are identical.
Matching points from the left and right codes are labeled with the same symbol.
For $N \le 5$ the points lie in the equatorial plane, and for $N \le 4$ there are repeated points.
The points lie on regular figures, except for $N=4$ where the points are
$\pm (1,0,0)$, $\pm \left( \frac{1}{\sqrt{3}} , \pm \sqrt{\frac{2}{3}} , 0 \right)$.

For $N=6$ the left and right codes consist of the 12 vertices of an icosahedron
(Fig.~\ref{fg12}).
Let these be the points
\beql{Eq20}
\la (0, 
\pm 1, \pm \tau ) ,~
\la ( \pm \tau , 0, \pm 1 ),~
\la ( \pm 1, \pm \tau , 0) ~,
\eeq
where $\la = 1 / \sqrt{ \tau +2}$.
The matching is obtained by mapping each point to its algebraic conjugate
(i.e. replacing $\sqrt{5}$ by $- \sqrt{5}$), and rescaling so the points
again lie on a unit sphere.
As already mentioned, the principal angles between each pair of these
planes are $\arcsin~ 1/\sqrt{5}$ and $\pi/2$, so $d_c^2 = 6/5$, $d_g^2 = 2.6824$.

\begin{figure}[htb]
\centerline{\psfig{file=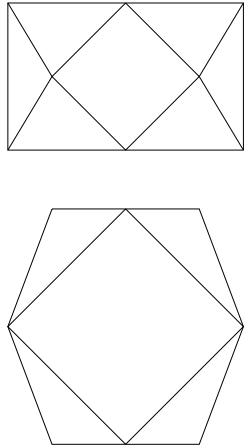,width=1.25in}}

\caption{Top and side views of anti-icosahedron.}
\label{fg13}
\end{figure}

There is a second set of six planes
with $d_c^2 = 6/5$, but with $d_g^2$ only equal to 2.0030.
Here the left and right codes form the vertices of what we shall call the ``anti-icosahedron'',
consisting of the points $\la( 0, \pm 1 , \pm \tau )$,
$\la ( \pm \tau , 0 , \pm 1)$,
$\la ( \pm \tau , \pm 1, 0)$ (compare \eqn{Eq20}), and shown in
Fig.~\ref{fg13}.
Topologically this
is a ``parallel bi-slit cuboctahedron'',
obtained by dividing two opposite square faces of a cuboctahedron into two triangles by parallel lines.
The matching sends
$\la(0, \pm 1, \tau )$ to $\la (0, \pm \tau, 1)$, $\la ( \tau , 0, \pm 1)$ to
$\la(1,0, \pm \tau )$ and
$\la (\pm \tau , 1, 0)$ to $\la ( \pm 1 , \tau , 0)$.

Four larger putatively $d_c$-optimal packings that the computer found
are also worth mentioning,
those of 10, 18, 48 and 50 planes.

For the 10-plane packing the left and right codes consist of the vertices of a decagonal
prism with coordinates
$$\pm \left( \sqrt{\frac{2}{3}} \cos ~r \theta , ~
\sqrt{\frac{2}{3}} \sin ~r \theta , ~ \sqrt{\frac{1}{3}} \right) ,~~
r=0, \dd, 9 ~,
$$
where $\theta = \pi /5$.
A typical point
$(\sqrt{2/3} \cos \, r\theta , ~ \sqrt{2/3} \sin\, r \theta ,~\sqrt{1/3} )$ in the left code is matched
with $(-1)^r$ $(\sqrt{2/3} \cos \, 3r \theta ,~ \sqrt{2/3} \sin 3r \theta , ~ \sqrt{1/3} )$ in the right code.
Then $d_c^2 = 10/9$, $d_g^2 = 1.5725$.
The most interesting property of this configuration is that, although three different sets of canonical
angles occur,
the chordal distance between every pair of planes is the same --- this is a regular simplex!

\begin{figure}[htb]
\centerline{\psfig{file=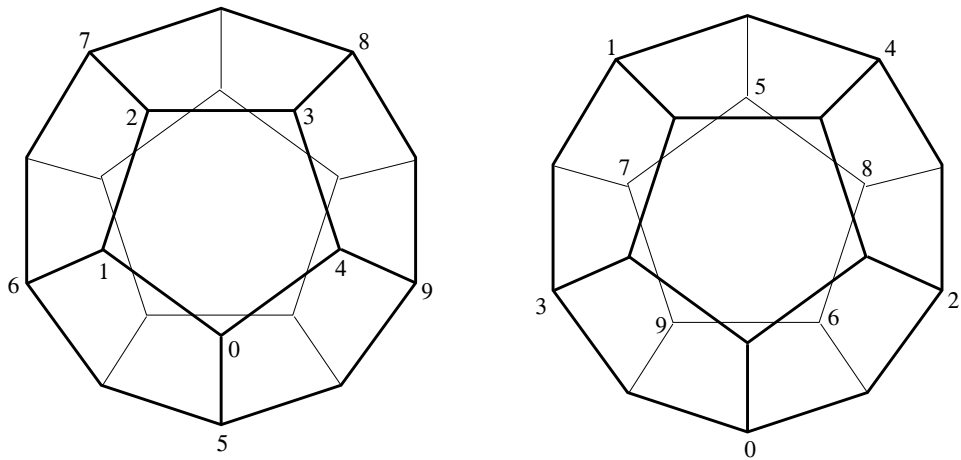,width=5in}}

\caption{Matching of vertices of dodecahedron.}
\label{fg14}
\end{figure}

We find it surprising that this is superior to the packing of ten planes obtained
by matching the vertices of a dodecahedron to their
algebraic conjugates, as shown in Figure~\ref{fg14}.
\paragraph{18 planes.}
Let $\sO$ denote the set of vertices of a regular octahedron.
The binocular code for the 18-plane packing is $\sO \times \sO := \{ ( \ell , r ) : \ell , r \in \sO \}$.
Projectively, the 18 two-spaces are the edges of a desmic triple\footnote{Two tetrahedra are called
{\em desmic} if they are in perspective from four distinct points --- these
form a third tetrahedron desmic with either.
The 18 edges of such a desmic triple are also the edges of a unique
other desmic triple.} of tetrahedra \cite{Ste79}, \cite{Hud05}, \cite{Cox50}.
We were pleased to see this configuration appear, since it was already familiar in the context
of quantum logic \cite{CK96}.
This set of planes can be constructed directly by taking the planes spanned by all pairs
of mutually perpendicular minimal vectors of the $D_4$ lattice (cf. \cite{SPLAG}).
In this form the planes have generator matrices such as
$$\left[
\matrix{ + & + & 0 & 0 \cr
0 & 0 & + & - \cr
}
\right] ~(12 ~{\rm planes} ) ,~~
\left[
\matrix{
+ & + & 0 & 0 \cr
+ & - & 0 & 0 \cr
}
\right] ~ (6 ~{\rm planes} ) ~.
$$
Three different sets of principal angles occur, $0$, $\pi /2$;
$\pi /4$, $\pi /4$;
and $\pi /2$, $\pi /2$; so that $d_c^2 =1$,
$d_g^2 = 1.2337$.
The automorphism group of this packing has structure
$[3,4,3].2$ and order 2304.
These 18 planes can also be described as the set of planes that meet the 24-cell
in one of its equatorial squares.
\paragraph{48 and 50 planes.}
Let $\sC$ denote the set of vertices of a cube.
Then the binocular codes for the 48- and 50-plane packings are respectively $\sO \times \sC \,\cup\, \sC \times \sO$ and $\sO \times \sO \,\cup \,\sC \times \sC$.
For both arrangements we have $d_c^2 = 2/3$, $d_g^2 = 0.7576$.
Both are believed to be optimal.
(We mention the 48-plane packing because of its greater symmetry.)
\paragraph{The matching problem.}
In order to maximize the minimal chordal distance between the planes,
the matching between the left and right codes should (from \eqn{Eq12})
be chosen so as to minimize the maximal value of $\cos \psi \cos \phi$.
Stated informally, the matching should be such that if two points are close together on one of the spheres then the points to which they are
matched should be far apart on the other sphere.
In the case of the icosahedron, for example (see Fig.~12), the matching sends adjacent vertices to non-adjacent vertices.
There is a unique way to do this.

At this point we were tempted to see if any new record packings
in $G(4,2)$ could be obtained by taking the $2N$ antipodal points corresponding
to a good $N$-line packing in $G(3,1)$, and matching them with themselves
in an optimal way.
However, it was not easy to see how to solve the matching problem.
Fortunately David Applegate (personal communication) found that it could be reformulated as an integer programming problem, as follows.

Given an antipodal set $S= \{P_1, \dd, P_{2N} \} \subseteq S^2$,
we wish to find a permutation $f$ of $S$ with the property that $f(-P_i) = f(P_i)$, for all $i$, and such that the minimal value of $1- P_i \cdot P_j f(P_i) \cdot f(P_j)$ $(i \neq j)$ is maximized.
If we do the maximization by binary search, we may define $\pi (P_i, P_j) =1$ if $f(P_i) = P_j$, or 0 otherwise, with the constraints
\begin{eqnarray*}
\pi (P_i , P_j) & = & \pi ( - P_i, - P_j ) ~, \\
\sum_{j=1}^{2N} \pi (P_i, P_j) & = & 1, ~~1 \le i \le 2N ~, \\
\sum_{i=1}^{2N} \pi (P_i, P_j) & = & 1 ,~~ 1 \le j \le 2N ~,
\end{eqnarray*}
and
$$
\pi (P_i, P_j) + \pi (P_k, P_\ell ) < 1 ~~\mbox{for all}~~
1 \le i,j,k, \ell \le 2N
~~\mbox{such that~ $1- P_i \cdot P_k ~ P_j \cdot P_\ell < M$} ~.
$$
There is a feasible solution if and only if $d_c^2 \ge M$.

Applegate kindly implemented this procedure, and used it to solve the matching
problem for our best packings in $G(3,1)$ and for various
polyhedra.
Unfortunately no new records have yet been obtained by this method.
\paragraph{Best packings for the geodesic metric.}
So far we have mostly discussed packings that attempt to maximize
the minimal chordal distance between planes.
We also computed packings for the geodesic distance,
and we shall now describe some of them.
In general, however, these are much less symmetric than the chordal-distance packings,
especially for more than 16 planes.

For $N=4,5$ the solutions are subsets of the 6-plane packing.
For $N=7$ the left code is a heptagonal anti-prism and the right code is an equatorial 14-gon.

For $N=12$ the binocular code is the set
$\{ \pm ( \ell , r ) \}$, where $\ell$ (resp. $r$) runs through
the vertices of a regular tetrahedron (resp. equatorial
equilateral triangle).

For $N=16$ the left and right codes are the union of two similarly oriented
square prisms.
For larger $N$ the geodesic packings do not seem so interesting.
For $N=18$, for example, the best packing has no non-trivial symmetries.
\paragraph{Hamiltonian paths.}
One further question remains to be discussed, that of arranging the planes in
a circuit in such a way that adjacent planes are close
together, for the ``Grand Tour'' application.
This turned out to be a much easier problem than finding the packings.
We handled it in two different ways.
For configurations such as the 48-plane arrangement, where there was an obvious notion of adjacency --- in this case, defining
two planes to be adjacent if the principal angles are $\pi /4$, $\pi /4$ --- we represent the packing by a graph with nodes representing the planes,
and look for a Hamiltonian cycle.
In less regular cases,
we convert the packing into a traveling salesman problem using chordal distance to define the distance between nodes,
and look for a minimal length circuit.
In both cases we were able to make use of
the travelling salesman programs of Applegate et~al. \cite{App95}, which
can handle 100-node graphs without difficulty.
Some of the packings in the database mentioned in Sect.~1 (those with
suffix {\tt .ham}) have been arranged in cycles in this way.
\begin{center}
~~ \\
~~ \\
\end{center}

The above reformulation in terms of binocular codes applies only to $G(4,2)$
(we realized from the beginning that this case would be special, since
$G(4,2)$ is the only Grassmannian space where the Riemannian metric is not unique
\cite{Lei61}).
In the next section we describe a second reformulation that applies to the
general case.

\section{Packing $n$-planes in $\RR^m$}
\hsp
Three observations contributed to the second reformulation.

(i)~We noticed (see Table~2) that for several examples of $N$-plane
packings in $G(4,2)$ the largest value of $d_c^2$ that we could attain was $N/(N-1)$, and that in every case this was an upper bound.
Further experimentation with other packings in $G(m,n)$ for $m \le 8$ led us
to guess an upper bound of
\beql{Eq30}
d_c^2 \le \frac{n(m-n)}{m} \cdot \frac{N}{N-1} ~,
\eeq
which again we could achieve for some small values of $N$.
The form of \eqn{Eq30} was suggestive of the Rankin bound for spherical codes \cite{SPLAG} or the Plotkin bound for binary codes \cite{MS77}.

(ii)~Investigation of the 18-plane chordal-distance packing revealed that,
with respect to chordal distance, this has the structure of a regular orthoplex (or generalized octahedron)
with 18 vertices.
Combining this with the fact that the 10-plane packing formed a regular
simplex, we had strong evidence that $G(4,2)$ should have an isometric
embedding (with respect to $d_c$) into $\RR^9$.

(iii)~A computer program was therefore written to determine the lowest
dimensions into which our library of packings
in $G(m,n)$ could be isometrically embedded.
More precisely, for a given set of $N$ planes in $G(m,n)$,
we searched for the smallest dimension $D$ such that
there are $N$ points in $\RR^D$
whose Euclidean distances coincide with the chordal
distances between our planes.

The results were a surprise:
it appeared that $G(m,n)$ with chordal distance could be isometrically embedded
into $\RR^D$, for $D= {\binom{m+1}{2}} -1$, independent of $n$.
Furthermore the points representing elements of $G(m,n)$ were observed
to lie on a sphere of radius $\sqrt{n(m-n)/2m}$ in $\RR^D$.

Aided by discussions with Colin Mallows, we soon found an explanation:
just associate to each $P \in G(m,n)$ the orthogonal projection map from
$\RR^m$ to $P$.
If $A$ is a generator matrix for $P$ whose rows are orthogonal unit vectors, then the projection is represented by the matrix
\beql{Eq29}
\sP = A^{tr} A~.
\eeq
$\sP$ is an $m \times m$ symmetric idempotent matrix, which is independent of the particular
orthonormal generator matrix used to define it.
Changing to a different coordinate frame in $\RR^m$ has the effect of conjugating
$\sP$ by an element of $O(m)$.
With the help of \eqn{Eq3}, we see that
\beql{Eq31}
{\rm trace} ~ \sP = n ~.
\eeq
Thus $\sP$ lies in a space of dimension
${\binom{m+1}{2}} -1$.

Let $\| ~ \|$ denote the $L_2$-norm of a matrix:
if $M= (M_{ij} )$, $1 \le i,j \le m$,
$$\| M \| = \sqrt{\sum_{i=1}^m \sum_{j=1}^m M_{ij}^2}
= \sqrt{{\rm trace} ~M^{tr} M} ~.
$$
For $P,Q \in G(m,n)$, with orthonormal generator matrices $A$, $B$,
and principal angles $\theta_1, \dd, \theta_n$, an elementary
calculation using \eqn{Eq3}, \eqn{Eq4} shows that
\begin{eqnarray}
\label{Eq32}
d_c^2 (P,Q) & = & n- ( \cos^2 \theta_1 + \cdots + \cos^2 \theta_n ) \nonumber \\
& = & n - {\rm trace} ~ A^{tr} A B^{tr} B \nonumber \\
& = & \frac{1}{2} \| \sP - \sQ \|^2 ~,
\end{eqnarray}
where $\sP$, $\sQ$ are the corresponding projection matrices.

Note that if we define the ``de-traced'' matrix
$\bar{\sP} = \sP - \frac{n}{m} I_m$, then ${\rm trace}~ \bar{\sP} =0$, and
$\| \bar{\sP} \|^2 = \frac{n(m-n)}{n}$.
We have thus established the following theorem.
\begin{theorem}\label{th2}
The representation of $n$-planes $P\in G(m,n)$ by their
projection matrices $\bar{\sP}$ gives an isometric
embedding of $G(m,n)$ into a sphere of radius $\sqrt{n(m-n)/n}$ in $\RR^D$,
$D= {\binom{m+1}{2}} -1$, with $d_c (P,Q) = \frac{1}{\sqrt{2}} \| \sP - \sQ \|$.
\end{theorem}

Thus chordal distance between planes is $1/ \sqrt{2}$ times the straight-line distance between the projection matrices (which explains our name for this metric).
The geodesic distance between the planes is $1/ \sqrt{2}$ times the geodesic distance between
the projection matrices, measured along the sphere in $\RR^D$.

Figure~\ref{fg15} attempts to display the embeddings of $G(m,0)$, $G(m,1) , \dd, G(m,m)$ in $\RR^{D+1}$.
Since $\| \sP - \frac{1}{2} I_m \|^2 = m/4$,
all the planes lie on the large sphere, centered at $\frac{1}{2} I_m$,
of radius $\sqrt{m} /2$.
The members of $G(m,n)$ lie on the intersection of the large sphere with the
plane ${\rm trace} (\sP) =n$, which intersection is itself a sphere in $\RR^D$ of
radius $\sqrt{n(m-n)/m}$ centered at $\frac{n}{m} I_m$.
A plane $P$ and its orthogonal complement $P^\perp$ are represented by antipodal points on the large sphere.
\begin{figure}[htb]
\begin{center}
\input Sglobe_text.tex
\end{center}

\vspace*{+.1in}
Figure~15:~\,Embedding of $G(m,0) , \dd, G(m,m)$ into large sphere of radius $\sqrt{m}/2$ in Euclidean space of dimension $m(m+1)/2$.  $G(m,n)$ lies on sphere of radius $\sqrt{n(m-n)/m}$ in $\RR^D$, $D= (m-1) (m+2) /2$.

\label{fg15}
\end{figure}

In contrast to this result, we briefly remark, without giving details,
that there is no way to embed $G(m,n)$ into Euclidean space of any dimension
so that the geodesic distance $d_g$ on $G(m,n)$ is represented by Euclidean distance in that space.
Of course the Pl\"{u}cker embedding, in which members of $G(m,n)$ are
represented by points in projective space of dimension
${\binom{m}{n}} -1$, also does not give a way to realize either $d_c$ or $d_g$ as Euclidean distance.
(Nor does the Nash embedding theorem \cite{Nas56}.)
Note also that the dimension of the Pl\"{u}cker embedding is in general much larger than the dimension of our embedding.

Since we have embedded $G(m,n)$ into a sphere of radius $\sqrt{n(m-n)/m}$ in $\RR^D$,
we can apply the Rankin bounds for spherical codes \cite{Ran55}, and deduce:
\paragraph{Corollary.}
{\em
(i)~The simplex bound:
for a packing of $N$ planes in $G(m,n)$,
\beql{Eq40}
d_c^2 \le \frac{n(m-n)}{m} \cdot \frac{N}{N-1} ~.
\eeq
Equality requires $N \le D+1 = {\binom{m+1}{2}}$, and occurs if and only if the $N$ points in $\RR^D$
corresponding to the planes form a regular `equatorial' simplex.

(ii)~The orthoplex bound: for $N > {\binom{m+1}{2}}$,
\beql{Eq41}
d_c^2 \le \frac{n(m-n)}{m} ~.
\eeq
Equality requires $N \le 2D = (m-1) (m+2)$, and occurs if the $N$ points form a
subset of the $2D$ vertices of a regular orthoplex.
If $N=2D$ this condition is also necessary.
}

\vspace*{+.1in}
Lemmons and Seidel (\cite{LS73b}, Theorem~3.6) give a bound for
equi-isoclinic packings in $G(m,n)$ which agrees with \eqn{Eq40};
of course our bound is more general.
The case $n=1$ of \eqn{Eq40} is given in Theorem~3 of
\cite{Ros96}.

We were happy to obtain confirmation of \eqn{Eq30}.
Since $d_c^2$ can never exceed $n$, we can also write
\beql{Eq40a}
d_c^2 \le \min \left\{
n, \frac{n(m-n)}{m} \frac{N}{N-1} \right\} ~.
\eeq

The Corollary allows us to establish the optimality of many of our packings.
In the range $m \le 16$, $n \le 3$, $N \le 55$,
there are over 750 cases which appear to meet \eqn{Eq40a} or
\eqn{Eq41}, in the sense that the ratio of $d_c^2$ to the bound is greater than .9999999~.
For $n=2$, for instance, the following are the packings of $N$ planes
in $G(m,2)$ for $m \le 10$ that appear to achieve
either \eqn{Eq40a} (listed before the semicolon) or \eqn{Eq41}
(after):
$$
\begin{array}{rl}
m=4: & 2-8, 10; 11-18 \\
m=5: & 4-11; 16, 17 \\
m=6: & 3-14; 22, 23 \\
m=7: & 6-18 ; 29 \\
m=8: & 4-21, 28; 37-44 \\
m=9: & 6-25; \\
m=10: & 5-28;
\end{array}
$$

However, this is not as meaningful as it at first seems.
Consider the packings of $N$ planes in $G(10,2)$.
The ratios of $d_c^2$ to the simplex bound for $N=21, \dd, 33$,
rounded to 14 decimal places, are as follows:
$$
\begin{array}{ll}
21&0.99999999999999 \\
22&1.00000000000000 \\
23&0.99999999999999 \\
24&0.99999999999999 \\
25&1.00000000000000 \\
26&0.99999999999999 \\
27&1.00000000000000 \\
28&0.99999999921354 \\
29&0.99999470290071 \\
30&0.99998661084736 \\
31&0.99961159256647 \\
32&0.99909699728979 \\
33&0.99854979246655
\end{array}
$$
From this it seems very likely that there exists a 27-plane packing meeting the simplex bound, but certainly further investigation is needed to determine if there
is a 28-plane simplex packing
in the neighborhood of the computer's approximate solution.
We have carried out such analyses in all the putative cases of equality for packing lines (see Sect.~6), for planes in $G(4,2)$ (as already discussed) and for 70 planes
in $G(8,4)$ (see below).
The other cases remain to be investigated.

In principle there is no difficulty in settling these questions.
All the chordal distances between the planes are specified.
So we could simply set up a (large) system of quartic equations with integer coefficients
for the entries in the generator matrices,
and ask if there is a real solution.

The chief difficulty in attempting to understand the computer-generated packings
is that there are usually very large numbers of equally good solutions.
The case of 18 planes in $G(4,2)$ is typical in this regard.
There is one exceptionally pleasing solution (described in Section~4), in which only
three different pairs of principal angles occur.
However, it seems that there is a roughly 29-dimensional manifold of solutions, in which a
typical solution, although still having $d_c^2 =1$, has an apparently random set of principal
angles.
(Our investigation of this question is not yet completed.)

The case of subspaces of dimension $n$ in $\RR^{2n}$ or $\RR^{2n+1}$ is especially interesting.
The largest possible arrangements of planes that could achieve the two bounds are:
\beql{Eq42}
\begin{array}{rrr@{~}lr@{~}l}
\multicolumn{1}{c}{m} & \multicolumn{1}{c}{n} &
\multicolumn{2}{c}{N({\rm simplex})} & \multicolumn{2}{c}{N({\rm orthoplex})} \\ [+.1in]
2 & 1 & 3 & \surd & 4 & \surd \\
3 & 1 & 6 & \surd & 10 & \\
4 & 2 & 10 & \surd & 18 & \surd \\
5 & 2 & 15  && 28 & \\
6 & 3 & 21 && 40 \\
7 & 3 & 28 & \surd & 54 & \\
8 & 4 & 36 & & 70 & \surd \\
\cdot & \cdot & \cdot
\end{array}
\eeq
Checks indicate that such a packing exists.
It is known that the orthoplex bound cannot be achieved by 10 planes
in $G(3,2)$,
while the other cases are undecided.
Our computer experiments strongly suggest that no set of 15 planes meets
the simplex bound in $G(5,2)$.

On the other hand it is possible to find
packings of 70 planes in $G(8,4)$ meeting the bound
\eqn{Eq41}.
With a considerable amount of effort (a long story, not told here), we determined several
examples, of which the following is the most symmetrical.
Let the coordinates be labeled $\infty , 0,1, \dd, 6$, and take two 4-planes generated by the vectors
\beql{Eq43}
\begin{array}{llll}
\{10000000, & 01000000, & 00100000, & 00001000 \} \\ [+.05in]
\{11000000, & 00101000, & 00010001, & 00000110 \}
\end{array} ~,
\eeq
respectively.
We obtain 70 planes from these by negating any even number of coordinates, and/or
applying the permutations $(0123456)$,
$(\infty 0) (16) (23) (45)$ and
$(124)(365)$.
The principal angles are $0,0, \frac{\pi}{2} , \frac{\pi}{2}$;
$\frac{\pi}{4}, \frac{\pi}{4} , \frac{\pi}{4} , \frac{\pi}{4}$;
or
$\frac{\pi}{2}, \frac{\pi}{2} , \frac{\pi}{2} , \frac{\pi}{2}$, so
$d_c^2 = 2$, $d_g^2 = \pi^2 /4$ (this is not even a local optimum with respect
to geodesic distance).
The full group, which is transitive on the planes,
is generated by the above operations and by the $8 \times 8$ Hadamard matrix
$$
\frac{1}{\sqrt{8}} \left[
\matrix{
--------\cr
-+--+-++\cr
---+-+++\cr
--+-+++-\cr
-+-+++--\cr
--+++--+\cr
-+++--+-\cr
-++--+-+\cr}
\right] ~,~
~~~
\begin{array}{c}
\mbox{($-$ means $-1$} \\ [+.05in]
\mbox{$+$ means $+1$)}
\end{array} ~.
$$
It has shape $2^8 \sA_8$ and order 5160960 ($\sA_n$ denotes an alternating
group of degree $n$).

A set of 28 planes in $G(7,3)$ meeting the simplex bound can be obtained as follows.
Label the coordinates $0, \dd, 6$, and let $v_r^\pm$ have components
1 at $2^r ( \bmod~7)$,
$\pm \sqrt{2}$ at $3.2^r ( \bmod~7)$, and 0's elsewhere,
for $r=0,1,2$.
Then the vectors $v_0^\pm$, $v_1^\pm$, $v_2^\pm$ (product of signs is even)
span four planes, and the full set of 28 is found by cycling the seven coordinates.\footnote{Note added April 1996.
P.~W. Shor and the third author have recently discovered that the packing
of 70 planes in $G(8,4)$ can be generalized to a packing of $m^2 +m-2$
planes in
$G(m,m/2)$ meeting the orthoplex bound,
whenever $m \ge 2$ is a power of 2;
and that the packing of 28 planes in $G(7,3)$ can be generalized to a packing of
$p(p+1)/2$ planes in $G(p, (p-1)/2)$ meeting the simplex bound, whenever
$p$ is a prime which is either 3 or congruent to $-1$ modulo 8.
These constructions will be described elsewhere.}

On the other hand, in spite of much effort, we have not been able to find a set of
40 planes in $G(6,3)$ that meet the bound.
If this were possible, we would obtain $d_c^2 =1.5$, whereas the best we have been able to
achieve is 1.49977.
We can get 1.5 with $N=34$ planes, and 1.49977 with $35 \le N \le 40$, but only 1.4297 with 41 planes, which suggests that 1.5 might be possible with 40 planes.
The discussion given earlier shows that this question could in principle be settled by seeing if a certain set of quadratic equations with integer coefficients has a real solution.
One formulation leads to 500 quartic equations in 360 unknowns, so this approach
seems hopeless at the present time.

It is not difficult to invent ways to construct sets of subspaces.
Several examples are given below.
However, many promising ideas have proved useless when confronted with the results
found by our program.
For this reason we have included two extensive tables of the best
chordal distances that we have found.
Table~\ref{ta3} gives $d_c^2$ for packings in $G(m,2)$, with $m \le 10$.
This will serve as a standard against which readers can test their own constructions.
$$
\begin{array}{c}
\hline
~~~~~~~~~~\mbox{Table~3 about here}~~~~~~~~~~ \\
\hline
\end{array}
$$

\paragraph{General constructions.}
The following are some promising general constructions.

(i)~A skew-symmetric conference matrix\footnote{For a list of the known orders of such matrices, see Table~7.1 of \cite{SeYa92}.} of order $4a$ yields a set of $4a$ unit vectors
in $\CC^{2a}$ with Hermitean inner products $\pm i / \sqrt{4a-1}$
(\cite{DGS}, Example~5.8), hence a set of $4a$ planes in $G(4a,2)$ with $d_c^2 = 4(2a-1) / (4a-1)$ that meet the simplex bound.

(ii)~Use the planes defined by the $n$-faces of a regular $m$-dimensional
polytope, or the Voronoi or Delaunay cells of a lattice, etc.
The initial results from this idea have been disappointing.
For example, the 96 two-dimensional faces of the 24-cell in $\RR^4$
define 16 different planes,
forming a packing with $d_c^2 = 8/9$, inferior to the best such packing (cf. Table~\ref{ta3}).

(iii)~Use the minimal vectors in a complex or quaternionic lattice to obtain packings
in $G(2a,2)$ or $G(4a,4)$.
For example, the 54 minimal vectors of the lattice $E_6^\ast$, regarded as a
three-dimensional lattice over the Eisenstein integers (\cite{SPLAG}, p.~127),
produce a packing of nine planes in $G(6,2)$ that meets the simplex bound.

(iv)~Restrict the search to generator matrices of 0's and 1's (as in \eqn{Eq43}) or
$+1$' and $-1$'s, or even to rows that are blocks in some combinatorial design, or
vectors in some error-correcting code (see the example in Section~6).

(v)~Let $C$ be an error-correcting code of length $m$ over $GF(2^n)$, for example
a Reed-Solomon code (cf. \cite{MS77}).
If we regard $GF(2^n)$ as a vector space of dimension $n$ over $GF(2)$, each codeword
yields an $n \times m$ matrix whose elements we may take to be $+1$'s and $-1$'s.
After discarding those of rank less than $n$, and weeding out duplicates, we
obtain a packing in $G(m,n)$.
The hexacode (\cite{SPLAG}, p.~82), for example, a code of length 6 over $GF(4)$
containing 64 codewords, produces 28 distinct planes in $G(6,2)$ with $d_c^2 = 3/4$.
Unfortunately Table~\ref{ta4} shows that the record is 1.2973.

(vi)~Choose a group $G$ with an $m$-dimensional representation, and a subgroup $H$ of index $N$ with an $n$-dimensional irreducible representation.
Find an $n$-dimensional subspace $V \subseteq \RR^m$ invariant under $H$, and take its orbit under $G$.
Many of the conference matrix constructions of line-packings described in
Section~6 are of this type (using $G= L_2 (q)$).

In some cases we have also used the optimizer to search for packings
with a specified group.
The following is a packing of 28 planes in $G(8,2)$ meeting the simplex bound that is an abstraction of a configuration found by the computer when searching for packings in $\RR^8$ invariant under the permutation $(0)(1,2,3,4,5,6,7)$.
Let $R: \CC^4 \to \RR^8$ map $(v_1, v_2, v_3, v_4)$ to
$({\rm Re} \, v_1, {\rm Im}\, v_1, {\rm Re} \, v_2, \dd, {\rm Im}\, v_4 )$, and let $\af = e^{2 \pi i/ 7}$.
The 28 planes are spanned by the following pairs of vectors:
$$
\begin{array}{lll}
R(0, \af^k , \af^{2k}, \af^{4k} ) , & R(0, i \af^k , i \af^{2k} , i \af^{4k} ) , & 0 \le k \le 6 ~, \\ [+.1in]
R(1,0, \af^{2k} , \af^{4k} ), & R(i, 0, i \af^{2k} , - i \af^{4k} ), & 0 \le k \le 6 ~, \\ [+.1in]
R(1, \af^k, 0, \af^{4k} ), & R(i, -i \af^k, 0, i \af^{4k} ), & 0 \le k \le 6~, \\ [+.1in]
R(1, \af^k , \af^{2k} , 0), & R(i, i \af^k , -i\af^{2k} , 0), & 0 \le k \le 6 ~.
\end{array}
$$
(The pattern of zeros and signs here suggests the Tetracode \cite{SPLAG}, p.~81.)
\paragraph{An application: apportioning randomness.}
The following is a potential application of packings in $G(m,n)$ for use in apportioning randomness, in the sense of producing large numbers of approximately random points from a few genuinely
random numbers.
We illustrate using the example of 70 planes in $G(8,4)$ described above.
Let $A_1, \ldots , A_{70}$ be generator matrices for them.
Suppose an exclusive resort wishes to distribute random points in $\RR^4$
to its 70 guests, perhaps for use as garage door keys.
Let $x= (x_1, \ldots, x_8)$ be a vector of eight independent Gaussian random
variates.
Then the hotel would assign $y_i = A_i x^{tr}$ to its $i$-th guest.
By maximizing the chordal distance between the planes we have minimized the correlation
between the $y_i$.
\section{Packing lines in higher dimensions}
\hsp
Table~\ref{ta4} shows the maximal angular separation found for packings of $N \le 50$ lines in $G(m,1)$ for $m \le 9$.

The following packings of $N$ lines in $G(m,1)$ for $m \le 10$ achieve
either the simplex or orthoplex bounds:
$m=3: ~ N=3,4,6,7$;
$m=4: ~ N=4,5,11,12$;
$m=5: ~ N=5,6,10,16$;
$m=6: ~ N=6,7,16,22$;
$m=7: ~ N=7,8,14,28$;
$m=8: ~ N=8,9$;
$m=9: ~ N=9,10,18,46-48$,
$m=10: ~ N=10,11,16$.
Most of these are described below.

As mentioned in Section~3, the best packing of 5 lines in $G(3,1)$ is a subset
of the best packing of 6 lines.
Table~4 shows similar phenomena in higher dimensions.
For example, the putatively best packing of 48 lines in $G(9,1)$ is so good that we cannot do better even when up through 8 lines are omitted from it.

In the rest of this section we discuss some of the entries in Table~\ref{ta4} that have the largest symmetry groups.
$$
\begin{array}{c}
\hline
~~~~~~~~~~\mbox{Table~4 about here}~~~~~~~~~~ \\
\hline
\end{array}
$$

\paragraph{Constructions from lattices.}
Let $P_1, \dd , P_a$ be a set of mutually touching spheres in a $d$-dimensional lattice packings.
If there are $2N$ further spheres in the packing, each of which touches all of
$P_1, \dd, P_a$, their centers form a $(d-a+1)$-dimensional antipodal spherical code with angular
separation ${\rm arcsec} (a+1)$ (\cite{SPLAG}, p.~340, Theorem~1).

In particular, the entries labeled ``1'' in Table~\ref{ta4} are obtained
from the centers of spheres that touch one sphere in the lattices
$D_4$, $D_5$, $E_6$ and $E_7$.
The entries labeled ``2'' are obtained from spheres that touch two spheres in the lattices
$D_4$, $D_5$, $E_6$, $E_7$, $E_8$, $\Lambda_9$ and the nonlattice packing
$P_{10b}$ (cf. \cite{SPLAG}, Chap.~1, Table~1.2).
We remark the set of 28 lines in $G(7,1)$ with angle ${\rm arccos} (1/3)$ obtained in this
way from $E_8$ is known to be unique:
see Chapter~14, Theorem~12 of \cite{SPLAG}.
This configuration of lines is derived in a different way in \cite{LS66}.

It is interesting to compare Table~\ref{ta4}
with the table of maximal sets of equiangular lines given in Lemmens and Seidel \cite{LS73} and
Seidel \cite{Sei94}.
Some arrangements appear in both tables,
for example the set of 28 lines in $G(7,1)$ just mentioned.
On the other hand the difference between the tables can be seen in dimension 8.
The maximal set of equiangular lines that exists in $G(8,1)$ has size 28,
with angle ${\rm arcsec} ~3$ (\cite{LS73}, Theorem~4.6).
However, Table~\ref{ta4} shows that there is a set of 28 lines (not
equiangular) in $G(8,1)$ with minimal angle arcsec~3.000511...~.
In view of the Lemmons-Seidel result (which uses the fact that in an $m$-dimensional
equiangular arrangement of $N$ lines with $N \ge 2m$ the secant must be an odd integer)
our set of 28 lines cannot be perturbed to give an equiangular set without decreasing the minimal angle.
\paragraph{From conference matrices.}
It follows from Theorem~6.3 of \cite{LS66} that if a symmetric conference matrix of order $q+1 \equiv 2$ $(\bmod~4)$ exists then there is an arrangement of $q+1$ equiangular
lines in $\RR^m$, $m= (q+1)/2$, with $d_c^2 = (q-1)/q$, meeting the simplex bound.
The corresponding entries are labeled ``3'' in the table.
Our program was able to find these packings for every prime power $q$ of this form below 100
except for 49 and 81.
\paragraph{From diplo-simplices.}
The entries labeled ``4'' are obtained by using all the vectors of shape
$\pm c (n^1, (-1)^n)$, where $c = 1/ \sqrt{n(n+1)}$.
In the notation of \cite{SPLAG}, Chap.~4, these are the minimal vectors in the translates
of the root lattice $A_n$ by the glue vectors $[1]$ and $[-1]$.
They form the vertices of a diplo-simplex \cite{CS91}, and have automorphism
group $2 \times \sA_{n+1}$.
These packings also meet the simplex bound.
\paragraph{From codes.}
Let $C$ be a binary code of length $m$, size $M$ and minimal distance $d$ which is
closed under complementation.
Writing the codewords as vectors of $\pm 1$'s, we obtain a packing of $M/2$ lines in
$G(m,1)$ with $d_c^2 = 4d(m-d)/m^2$.
For example, a shortened Hamming code of length 10 with $M=32$, $d=4$ gives a packing of 16 lines in $G(10,1)$ meeting the simplex bound.
This construction provides a rich supply of good packings.
The Nordstrom-Robinson code \cite{FST}, \cite{HKCSS}, for example, yields a packing of 128 lines in $G(16,1)$ with $d_c^2 = 15/16$.
\paragraph{40 lines in $G(4,1)$.}
The 80 antipodal points are $(\xi^{\mu+1/2}, 0)$, $(0, \xi^{\nu+1/2} )$,
$(a \xi^{2 \mu} , b \xi^{2 \nu +1} )$,
$(a \xi^{2 \mu+1} , b \xi^{2 \nu} )$,
$(b \xi^{2 \mu} , a \xi^{2 \nu} )$, $(b \xi^{2\mu +1} , a \xi^{2 \nu+1} )$,
where $\xi = e^{2 \pi i /8}$, $\mu, \nu =0, \dd, 3$, $a=2^{-1/4}$,
$b= \sqrt{1-a^2}$.
The group of these 80 points is the group $\frac{1}{2} ( \sD_{16} + \sD_{16} ) \cdot 2$ of order 256 (see \cite{CS96}).
Here $\sD_m$ denotes a dihedral group of order $m$.
\paragraph{22 lines in $G(6,1)$.}
We use the vertices of a hemi-cube (an alternative way to describe the 16-line packing)
together with the six coordinate axes.
\paragraph{63 and 64 lines in $G(7,1)$.}
The putatively best packing of 63 lines in $G(7,1)$ (just beyond the range of Table~\ref{ta4}) has angular separation $60^\circ$, and is formed from the 126 minimal vectors of the lattice $E_7$.
The automorphism group of this set of 126 points is the Weyl group
$\sW (E_7)$, of order $2^{10} .3^4 . 5.7$ (cf. \cite{SPLAG}, Chap.~4, \S8.3).

The best packing found of 64 lines in $G(7,1)$ also has an unusually large group, of order $2 | \sW (E_6) | = 2^8 . 3^4 . 5$.
This packing can be obtained as follows.
The largest gaps between the packing of 63 lines occur
in the directions of the minimal vectors of the dual lattice
$E_7^\ast$.
Let $u$ be such a vector, with $u= \slashedfrac{1}{2}\, v_6$, $v_6 \in E_7$,
$v_6 \cdot v_6 = 6$.
Taking a coordinate frame in which the last coordinate is in the $v_6$ direction,
we obtain the minimal vectors of $E_7$ in the form $(w,0)$, $w \in E_6$, $w \cdot w = 2$ (72 vectors), and
$(x, \pm \sqrt{2/3} )$, $x \in E_6^\ast$, $x \cdot x = 4/3$
(54 vectors).
We now adjoin $\pm \sqrt{6}$ and rescale by multiplying the last coordinate by
$\sqrt{3/4}$, obtaining our final configuration of vectors
$(0, \pm 3 / \sqrt{2})$, $(w,0)$, $(x, \pm 1/ \sqrt{2} )$, which is isomorphic to
what the computer found.
The rescaling has compressed the $E_7$ lattice in the $v_6$ direction until the angle
between $v_6$ and the $(x, 1 /\sqrt{2} )$ vectors is the same as the minimal angle between the $(x, 1/ \sqrt{2} )$ and $(w,0)$ vectors.
This angle,
the minimal angle in the packing, is ${\rm arccos} ( \sqrt{3/11}) = 58.5178^\circ$.
\paragraph{36 lines in $G(8,1)$.}
The 72 antipodal points consist of all permutations of the two vectors
$c(7^2, -2^7)$
and $c(-7^2, 2^6)$, where $c = 1 / \sqrt{126}$.
These are the minimal vectors in the translates of the $A_8$ lattice by the glue
vectors $[2]$ and $[-2]$.

The configuration of 16 lines in $G(5,1)$ is similarly obtained from the translates
of $A_5$ by $[1]$, $[3]$ and $[5]$.
\paragraph{39 lines in $G(12,1)$.}
Take the 13 lines of a projective plane of order 3, described by vectors
of four 1's and nine 0's, and replace two of the 1's by $-1$'s in all possible ways,
obtaining 78 vectors in $\RR^{12}$ with angle arcsec~4.
The group is $2 \times L_3 (3)$, of order $2^5 3^3 13$.

If instead we replace any {\em odd} number of the 1's by $-1$'s,
we obtain a putatively optimal packing of 52 lines in $G(13,1)$, with angle
arcsec 4 and group $2 \times L_3 (3).2$ of order $2^6 3^3 13$.

As in all these examples, we were amazed that the program was able to
discover such beautiful configurations.
\clearpage
\subsection*{Acknowledgements}
\hsp
We are grateful to Daniel Asimov, Andreas Buja, Dianne Cook, Simon Kochen, Warren Smith, and especially
Colin Mallows for helpful discussions, and to David Applegate for help in solving the
matching problem described in Section~4.
\begin{table}[htb]
\caption{Best packings found of $N \le 50$ planes in $G(m,2)$, $m \le 10$.
The entry gives $d_c^2$.}

$$\mbox{\small $\begin{array}{cccccccc}
N \setminus m & 4 & 5 & 6 & 7 & 8 & 9 & 10 \\ [+.1in]
~3 & 1.5000 & 1.7500 & 2.0000 & 2.0000 & 2.0000 & 2.0000 & 2.0000 \\
~4 & 1.3333 & 1.6000 & 1.7778 & 1.8889 & 2.0000 & 2.0000 & 2.0000 \\
~5 & 1.2500 & 1.5000 & 1.6667 & 1.7854 & 1.8750 & 1.9375 & 2.0000 \\
~6 & 1.2000 & 1.4400 & 1.6000 & 1.7143 & 1.8000 & 1.8667 & 1.9200 \\
~7 & 1.1667 & 1.4000 & 1.5556 & 1.6667 & 1.7500 & 1.8148 & 1.8667 \\
~8 & 1.1429 & 1.3714 & 1.5238 & 1.6327 & 1.7143 & 1.7778 & 1.8286 \\
~9 & 1.1231 & 1.3500 & 1.5000 & 1.6071 & 1.6875 & 1.7500 & 1.8000 \\
10 & 1.1111 & 1.3333 & 1.4815 & 1.5873 & 1.6667 & 1.7284 & 1.7778 \\
11 & 1.0000 & 1.3200 & 1.4667 & 1.5714 & 1.6500 & 1.7111 & 1.7600 \\
12 & 1.0000 & 1.3064 & 1.4545 & 1.5584 & 1.6364 & 1.6970 & 1.7455 \\
13 & 1.0000 & 1.2942 & 1.4444 & 1.5476 & 1.6250 & 1.6852 & 1.7333 \\
14 & 1.0000 & 1.2790 & 1.4359 & 1.5385 & 1.6154 & 1.6752 & 1.7231 \\
15 & 1.0000 & 1.2707 & 1.4286 & 1.5306 & 1.6071 & 1.6667 & 1.7143 \\
16 & 1.0000 & 1.2000 & 1.4210 & 1.5238 & 1.6000 & 1.6593 & 1.7067 \\
17 & 1.0000 & 1.2000 & 1.4127 & 1.5179 & 1.5937 & 1.6528 & 1.7000 \\
18 & 1.0000 & 1.1909 & 1.4048 & 1.5126 & 1.5882 & 1.6471 & 1.6941 \\
19 & 0.9091 & 1.1761 & 1.3948 & 1.5078 & 1.5833 & 1.6420 & 1.6889 \\
20 & 0.9091 & 1.1619 & 1.3888 & 1.5026 & 1.5789 & 1.6374 & 1.6842 \\
21 & 0.8684 & 1.1543 & 1.3821 & 1.4987 & 1.5750 & 1.6333 & 1.6800 \\
22 & 0.8629 & 1.1419 & 1.3333 & 1.4912 & 1.5714 & 1.6296 & 1.6762 \\
23 & 0.8451 & 1.1332 & 1.3333 & 1.4859 & 1.5680 & 1.6263 & 1.6727 \\
24 & 0.8372 & 1.1251 & 1.3326 & 1.4790 & 1.5638 & 1.6232 & 1.6696 \\
25 & 0.8275 & 1.1178 & 1.3229 & 1.4725 & 1.5594 & 1.6204 & 1.6667 \\
26 & 0.8144 & 1.1113 & 1.3151 & 1.4666 & 1.5556 & 1.6177 & 1.6640 \\
27 & 0.8056 & 1.1045 & 1.3071 & 1.4606 & 1.5556 & 1.6154 & 1.6615 \\
28 & 0.8005 & 1.0989 & 1.2987 & 1.4531 & 1.5556 & 1.6118 & 1.6593 \\
29 & 0.7889 & 1.0937 & 1.2887 & 1.4286 & 1.5455 & 1.6083 & 1.6571 \\
30 & 0.7809 & 1.0875 & 1.2804 & 1.4234 & 1.5398 & 1.6049 & 1.6552 \\
31 & 0.7760 & 1.0822 & 1.2675 & 1.4167 & 1.5342 & 1.6011 & 1.6527 \\
32 & 0.7691 & 1.0766 & 1.2588 & 1.4106 & 1.5304 & 1.5978 & 1.6501 \\
33 & 0.7592 & 1.0722 & 1.2526 & 1.4038 & 1.5244 & 1.5935 & 1.6476 \\
34 & 0.7549 & 1.0671 & 1.2447 & 1.3978 & 1.5216 & 1.5911 & 1.6448 \\
35 & 0.7489 & 1.0640 & 1.2430 & 1.3915 & 1.5158 & 1.5893 & 1.6431 \\
36 & 0.7477 & 1.0596 & 1.2345 & 1.3843 & 1.5086 & 1.5885 & 1.6414 \\
37 & 0.7286 & 1.0519 & 1.2272 & 1.3784 & 1.5000 & 1.5816 & 1.6364 \\
38 & 0.7198 & 1.0462 & 1.2239 & 1.3722 & 1.5000 & 1.5768 & 1.6334 \\
39 & 0.7095 & 1.0404 & 1.2206 & 1.3663 & 1.5000 & 1.5730 & 1.6305 \\
40 & 0.7066 & 1.0355 & 1.2149 & 1.3606 & 1.5000 & 1.5695 & 1.6289 \\
41 & 0.6992 & 1.0302 & 1.2115 & 1.3575 & 1.5000 & 1.5665 & 1.6242 \\
42 & 0.6948 & 1.0256 & 1.2079 & 1.3550 & 1.5000 & 1.5639 & 1.6213 \\
43 & 0.6844 & 1.0207 & 1.2037 & 1.3449 & 1.5000 & 1.5616 & 1.6192 \\
44 & 0.6831 & 1.0159 & 1.2007 & 1.3407 & 1.5000 & 1.5588 & 1.6174 \\
45 & 0.6809 & 1.0122 & 1.1969 & 1.3354 & 1.4839 & 1.5570 & 1.6160 \\
46 & 0.6793 & 1.0078 & 1.1941 & 1.3321 & 1.4821 & 1.5480 & 1.6099 \\
47 & 0.6732 & 1.0042 & 1.1924 & 1.3265 & 1.4574 & 1.5412 & 1.6072 \\
48 & 0.6667 & 1.0001 & 1.1907 & 1.3232 & 1.4490 & 1.5353 & 1.6045 \\
49 & 0.6667 & 0.9960 & 1.1873 & 1.3209 & 1.4430 & 1.5304 & 1.6019 \\
50 & 0.6667 & 0.9910 & 1.1841 & 1.3150 & 1.4370 & 1.5257 & 1.6003
\end{array}$}
$$
\label{ta3}
\end{table}
\begin{table}[htb]
\caption{Maximal angular separation found for $N \le 50$ lines in $G(m,1)$, $m \le 9$.}

$$\mbox{\small $\begin{array}{clllllll}
N \setminus m & \multicolumn{1}{c}{3} & \multicolumn{1}{c}{4} & \multicolumn{1}{c}{5} & \multicolumn{1}{c}{6} & \multicolumn{1}{c}{7} & \multicolumn{1}{c}{8} & \multicolumn{1}{c}{9} \\ [+.1in]
~3 & 90.0000 & 90.0000 & 90.0000 & 90.0000 & 90.0000 & 90.0000 & 90.0000 \\
~4 & 70.5288^2 & 90.0000 & 90.0000 & 90.0000 & 90.0000 & 90.0000 & 90.0000 \\
~5 & 63.4349^6 & 75.5225^4 & 90.0000 & 90.0000 & 90.0000 & 90.0000 & 90.0000 \\
~6 & 63.4349^3 & 70.5288^2 & 78.4630^4 & 90.0000 & 90.0000 & 90.0000 & 90.0000 \\
~7 & 54.7356^6 & 67.0213 & 73.3689 & 80.4059^4 & 90.0000 & 90.0000 & 90.0000 \\
~8 & 49.6399^6 & 65.5302 & 70.8039 & 76.0578 & 81.7868^4 & 90.0000 & 90.0000 \\
~9 & 47.9821^6 & 64.2619 & 70.5288 & 73.8437 & 78.4630 & 82.8192^4 & 90.0000 \\
10 & 46.6746^6 & 64.2619 & 70.5288^2 & 73.6935 & 76.3454 & 79.4704 & 83.6206^4 \\
11 & 44.4031 & 60.0000 & 67.2543 & 71.5651 & 75.0179 & 77.8695 & 80.6204 \\
12 & 41.8820^6 & 60.0000^1 & 67.0213 & 71.5651 & 74.1734 & 76.6050 & 79.4704 \\
13 & 39.8131 & 55.4646 & 65.7319 & 70.5288 & 73.8979 & 76.1645 & 77.9422 \\
14 & 38.6824 & 53.8376 & 65.7241 & 70.5288 & 73.8979^3 & 75.0349 & 77.2382\\
15 & 38.1349^6 & 52.5016 & 65.5302 & 70.5288 & 71.5678 & 74.3318 & 76.5006\\
16 & 37.3774^6 & 51.8273 & 63.4349^5 & 70.5288^2 & 70.9861 & 74.1005 & 75.9638\\
17 & 35.2353 & 50.8870 & 61.2551 & 68.1088 & 70.5926 & 73.1371 & 75.9638\\
18 & 34.4088 & 50.4577 & 61.0531 & 67.3744 & 70.5527 & 72.7464 & 75.9638^3\\
19 & 33.2115 & 49.7106 & 60.0000 & 67.3700 & 70.5288 & 72.0756 & 74.4577\\
20 & 32.7071 & 49.2329 & 60.0000^6 & 67.0996 & 70.5288 & 71.6706 & 74.2278\\
21 & 32.2161 & 48.5479 & 57.2025 & 67.0213 & 70.5288 & 71.3521 & 73.7518\\
22 & 31.8963 & 47.7596 & 56.3558 & 65.9052^5 & 70.5288 & 71.0983 & 73.1894\\
23 & 30.5062 & 46.5104 & 55.5881 & 63.6744 & 70.5288 & 70.7720 & 72.7488\\
24 & 30.1628 & 46.0478 & 55.2279 & 63.6122 & 70.5288 & 70.6027 & 72.6547\\
25 & 29.2486 & 44.9471 & 54.8891 & 62.4240 & 70.5288 & 70.5490 & 72.3124\\
26 & 28.7126 & 44.3536 & 54.2116 & 61.7377 & 70.5288 & 70.5432 & 72.1763\\
27 & 28.2495 & 43.5530 & 53.5402 & 61.4053 & 70.5288 & 70.5392 & 71.6650\\
28 & 27.8473 & 43.1566 & 53.2602 & 60.5276 & 70.5288^2 & 70.5322 & 71.5794\\
29 & 27.5244 & 42.6675 & 53.0180 & 60.1344 & 66.7780 & 70.5288 & 71.5175\\
30 & 26.9983 & 42.2651 & 52.7812 & 60.0213 & 65.7563 & 70.5288 & 71.5175\\
31 & 26.4987 & 42.0188 & 52.4120 & 60.0000 & 65.1991 & 70.5288 & 70.8508\\
32 & 25.9497 & 41.9554 & 52.3389 & 60.0000 & 64.7219 & 70.5288 & 70.7437\\
33 & 25.5748 & 41.4577 & 52.2465 & 60.0000 & 64.6231 & 69.3203 & 70.6940\\
34 & 25.2567 & 40.9427 & 51.8537 & 60.0000 & 64.6231 & 69.1688 & 70.6512\\
35 & 24.8702 & 40.7337 & 51.8273 & 60.0000 & 64.6231 & 69.0752 & 70.6337\\
36 & 24.5758 & 40.6325 & 51.8273 & 60.0000^1 & 64.6231 & 69.0752^5 & 70.5864\\
37 & 24.2859 & 40.4486 & 51.8273 & 57.6885 & 62.3797 & 67.7827 & 70.5695\\
38 & 24.0886 & 40.4419 & 50.3677 & 57.1057 & 62.1435 & 67.3835 & 70.5571\\
39 & 23.8433 & 39.6797 & 50.0611 & 56.8357 & 61.7057 & 67.1387 & 70.5443\\
40 & 23.3293 & 39.0236^5 & 49.5978 & 56.0495 & 61.3792 & 66.3815 & 70.5288\\
41 & 22.9915 & 38.5346 & 49.2600 & 55.8202 & 61.1630 & 65.9282 & 70.5288\\
42 & 22.7075 & 38.3094 & 48.6946 & 55.6160 & 60.8232 & 65.8166 & 70.5288\\
43 & 22.5383 & 37.7833 & 48.4030 & 55.3981 & 60.5193 & 65.2885 & 70.5288\\
44 & 22.2012 & 37.3474 & 48.0955 & 55.1259 & 60.3623 & 65.2422 & 70.5288\\
45 & 22.0481 & 37.1198 & 47.7723 & 54.9980 & 60.2282 & 64.7476 & 70.5288\\
46 & 21.8426 & 36.9997 & 47.3753 & 54.9858 & 60.1101 & 64.4007 & 70.5288\\
47 & 21.6609 & 36.5952 & 47.0323 & 54.7356 & 60.0433 & 64.1542 & 70.5288\\
48 & 21.4663 & 36.3585 & 46.7105 & 54.5940 & 60.0116 & 63.8846 & 70.5288^2\\
49 & 21.1610 & 36.1369 & 46.4345 & 54.5031 & 60.0000 & 63.4849 & 68.0498\\
50 & 20.8922 & 36.0754 & 46.1609 & 54.3191 & 60.0000^1 & 63.1527 & 67.7426
\end{array}$}
$$

{\bf Key.}
1,2: from sphere packings;
3: from conference matrices;
4: diplo-simplex;
5: described below;
6: described in Section~3.
\label{ta4}
\end{table}

\end{document}